\pdfoutput=1
\documentclass[opre,nonblindrev]{informs_modified} 

\DoubleSpacedXI 

\usepackage{endnotes}
\let\footnote=\endnote

\newcommand{\bN}{\mathbb{N}}
\newcommand{\bE}{\mathbf{E}}
\newcommand{\bP}{\mathbf{P}}
\newcommand{\Var}{\text{Var}}

\newcommand{\ratiok}{1 - \sqrt{\frac{2}{\pi}} \frac{1}{\sqrt{k}} + O(\frac{1}{k})}

\usepackage{natbib}
 \bibpunct[, ]{(}{)}{,}{a}{}{,}%
 \def\bibfont{\small}%
\TheoremsNumberedThrough     
\ECRepeatTheorems

\EquationsNumberedThrough    

\begin{document}

\RUNAUTHOR{Wang, Truong, and Bank}
\RUNTITLE{Online Advance Admission Scheduling}

\TITLE{Online Advance Admission Scheduling for Services with Customer Preferences}
\ARTICLEAUTHORS{%
\AUTHOR{Xinshang Wang, Van-Anh Truong} \AFF{Department
of Industrial Engineering and Operations Research, Columbia
University, New York, NY, USA, \EMAIL{xw2230@columbia.edu, vatruong@ieor.columbia.edu}
\URL{}}
\AUTHOR{David Bank, MD, MBA}
\AFF{Department of Pediatrics, NYPH \/ Morgan Stanley Children's Hospital, Columbia University Medical Center, New York, NY, USA, \EMAIL{deb40@columbia.edu}} 
}

\ABSTRACT{We study web and mobile applications that are used to schedule advance service, from medical appointments to restaurant reservations.  We model them as online weighted bipartite matching problems with non-stationary arrivals. We propose new algorithms with performance guarantees for this class of problems.  Specifically, we show that the expected performance of our algorithms is bounded below by $\ratiok$ times that of an optimal offline algorithm, which knows all future information upfront, where $k$ is the minimum capacity of a resource.  This is the tightest known lower bound. 
This performance analysis holds for any Poisson arrival process. Our algorithms can also be applied to a number of related problems, including display ad allocation problems and revenue management problems for opaque products. 
We test the empirical performance of our algorithms against several well-known heuristics by using appointment scheduling data from a major academic hospital system in New York City. The results show that the algorithms exhibit the best performance among all the tested policies. In particular, our algorithms are $21\%$ more effective than the actual scheduling strategy used in the hospital system according to our performance metric.}

\maketitle
\section{Introduction}
We study advance admission scheduling decisions in service systems.  Advance admission scheduling decisions are those that determine specific times for customers' arrival to a facility for service.   Advance admission scheduling is used in many service industries. Restaurants reserve tables for customers who call in advance. Healthcare facilities reserve appointment slots for patients who request them. Airlines reserve flight seats for those who purchase flight tickets. Advance admission scheduling enables service providers to better match capacity with demand because they control customers' actual arrivals to service facilities.  

We formulate and analyze a model that generally captures such admission scheduling systems. For concreteness, we focus on the example of MyChart, a digital admission scheduling application developed by Epic System. Epic is an electronic medical records company that is managing the records of millions of health care providers and more than half of the patient population in the U.S. \citep{epicMyChart}. Epic deploys MyChart to perform online scheduling of appointments through internet portals. The use of applications like MyChart is part of a general trend in healthcare towards providing electronic access to service through web and mobile applications 
\citep{technologyAdvice2015}.   

When a patient schedules an appointment over a web portal, MyChart first asks the patient for the type of visit desired, whether it is for a physical exam, a consultation, a flu shot, etc.  Next, it asks for the beginning and end of the range of preferred dates.  It then shows a menu with a check box for morning and afternoon session for each day in the preferred date range.  Patients can select one or more preferred sessions.  Finally, MyChart either offers the patient one or more appointments, or states that no appointment can be found. We can conceive of many variations over this basic interface.




Consider the following model of advance admission scheduling that captures MyChart as an example. There are multiple service providers.  Each provider offers a number of service sessions over a continuous, finite horizon.  We call a session associated with a single provider a \emph{resource.} Let $n$ be the number of resources available over the horizon. Each resource $j$ can serve $C_j$ customers. We call $C_j$ the \emph{capacity} of resource $j$.  
Each resource $j$ must be booked by time $t_j$ or it perishes at time $t_j$. There are $m$ customer \emph{types}. Patients of type $i$, $i=1,\ldots,m$, arrive according to some known non-homogeneous Poisson process and make reservations through any of the modes made available by the provider, web, phone, or mobile. 
A patient of type $i$ generates a \emph{reward} of $r_{ij}$ when served with a unit of resource $j$.  
We assume that the type of customers can be observed at the time that they arrive to make an appointment, through the pattern of preferences that they indicate and any data stored in the system on their profiles. 
We require that customers arriving at time $t$ have weight $0$ for all resources $j$ that perish at time $t_j < t$.  The number of customer types can be kept finite by discretizing the horizon but this number can be very large.  We will discuss this point shortly.  When a customer arrives, a unit of an available resource must be assigned to her, or she must be rejected.  Each unit of a resource can be assigned to at most one customer. We allow no-shows and the practice of overbooking to compensate for the effect of no-shows. The objective of the problem is to allocate the resources to the customers to maximize the expected total reward of the allocation.

Our advance reservation model is essentially an online weighted bipartite matching problem. The resources in our model, when partitioned into units, can be seen as nodes on one side of a bipartite graph. All the customers correspond to nodes on the other side that are arriving online. The type of each arriving customer is determined by a time-varying distribution. 

This resource allocation model can be found in many other applications.  We summarize three such applications below.
\begin{description}
\item[\bf Ad allocation.] In a typical display ad allocation problem, e-commerce companies aim at tailoring display ads for each type of customers. Each ad, which corresponds to a resource, is often associated with a maximum number of times to be displayed. Knowing the arrival rates of future customers, the task is to make the most effective matching between ads and customers. 

\item[\bf Single-leg revenue  management.] A special case of our model is the classic single-leg revenue management problem in which all resources to be allocated are available at the same time. Customers who bring a higher reward correspond to higher-fare classes. The decision is how to admit or reject customers, given the time remaining until the flight and the current inventory of available seats. 

\item[\bf Management of opaque products.] Internet retailers such as Hotwire or Priceline often offer a buyer an under-specified or \emph{opaque product}, such as a flight ticket, with certain details such as the exact flight timing or the name of the airline withheld until after purchase. We assume that demand for each opaque product is exogenous and independent of the availability of other products. When demand occurs, a decision is made to assign a specific product to that demand unit.  Knowing the arrival rates of all demands, we want to maximize the total expected revenue by strategically assigning specific products. 
\end{description}

Our contributions in this work are as follows:

\begin{itemize}
\item We provide a general, high-fidelity model of advance admission scheduling that captures customer preferences across different resources. 
We allow non-stationary arrivals and no-shows. We model the advance admission scheduling problem as an online weighted bipartite matching problem with non-stationary arrivals and propose new algorithms with guarantees on the relative performance.

\item We prove the tightest known performance bound for the online matching problem with non-stationary stochastic arrivals.  Specifically, we prove that a primitive algorithm, which we call the \emph{Separation Algorithm}, has expected performance that is bounded below by $\ratiok$ times that of an optimal offline algorithm, which knows all future information upfront, where $k$ is the minimum capacity of any resource. 
Our performance bound improves upon the lower bound of Alaei,
Hajiaghayi and Liaghat (2012). Moreover, it is close to an upper bound on the performance of the Separation Algorithm that the same authors found.

We obtain our bound by analyzing a novel \emph{bounded Poisson process}.  This is a Poisson process to which we apply a sequence of reflecting barriers.  The process arises in the dual of an optimization problem that characterizes our performance bound.  The behavior of this process is very complex, with no known closed-form description.  We managed to obtain a closed-form approximate characterization of the process.

\item  We improve on the Separation Algorithm by devising a novel bid-price-based algorithm, called the \emph{Marginal Allocation Algorithm}, that is much more practical.  First, the Marginal Allocation Algorithm is non-randomized, therefore more stable.  Second, it is fair in the sense that it never rejects a high-priority customer but accepts a low-priority customer, assuming that their arrival times and preferences are the same.  We prove that the Marginal Allocation Algorithm\ has the same theoretical performance guarantee as the Separation Algorithm.  In addition, in numerical experiments, we show that it achieves much better practical performance.


\item We test the empirical performance of our algorithm against several well-known heuristics by using appointment scheduling data from a department within a major academic hospital system in New York City. The results show that our scheduling algorithms perform the best among all tested policies. In particular, our algorithm is $21\%$ more effective than the actual scheduling strategy used in the hospital system according to our performance metric.
\end{itemize}


\section{Literature Review}

\subsection{Appointment Scheduling}

Our work is related to the literature on appointment scheduling. This area has been studied intensively in recent years \citep{guerriero2011operational,may2010surgical,cardoen2010operating,gupta2007surgical}. A large part of this literature considers {\it intra-day} scheduling, in which the number of patients to be treated on each day is given or is exogenous, and the task is to determine an efficient sequence of start times for their appointments. Another part of the literature considers {\it multi-day} scheduling, in which patients are dynamically allocated to appointment days. Some works in this literature focus on the number of patients to be served today, with the rest of the patients remaining on a waitlist until the next day. This paradigm is called {\it allocation scheduling}. See, for example, \citet{Huh2013}, \citet{Min20101091,Huh2010}, \citet{Gerchak1996}. Recently, more works have focused on the problem of directly scheduling patients into future days. This paradigm is called {\it advance scheduling}. This paper considers an advance scheduling model with multiple patient classes. 
In the literature of advance scheduling, \citet{truong2014optimal} first studies the analytical properties of a two-class advance scheduling model and gives efficient solutions to an optimal scheduling policy. For the multi-class model, no analytical result is known so far. \citet{Gocgun20122323} and \citet{Patrick2008} propose heuristics based on approximate dynamic programming for these problems, but have not characterized the worst-case performance of these heuristics. We propose online scheduling policies with performance guarantees for a very general multi-class advance scheduling problem.

Our advance scheduling model captures the preferences of patients in a general way.  Patient preferences are an important consideration in most out-patient scheduling systems. In the literature considering patient preferences, \citet{doi:10.1287/opre.1080.0542} consider a single-day scheduling model where each arriving patient picks a single slot with a particular physician, and the clinic accepts or rejects the request. 
Our model can be seen as a multi-period generalization of their work. We also characterize the theoretical performance in an online setting, whereas they use stochastic dynamic programming as the modeling framework and develop heuristics.
\citet{doi:10.1287/opre.2014.1286}  study how to offer sets of open appointment slots to a stream of arriving patients over a finite horizon of multiple days, given that patients have preferences for slots that can be captured by the multinomial logit model.  Their work is strongly influenced by assortment planning problems.  An important observation, which was first made by \cite{doi:10.1287/opre.1080.0542}, is that there is a fundamental difference between many advance admission scheduling problems and assortment planning problems.  In admission scheduling, we can often work with revealed preferences, whereas in assortment planning problems, decisions are made with knowledge only of a distribution of customer preferences. Working with revealed preferences allows for a more efficient allocation of service compared to working with opaque preferences.  It also leads to more analytically tractable models.


\subsection{Online Resource Allocation}

Our work is closely related to works on online matching problems. Traditionally, the online bipartite matching problem studied by \citet{Karp:1990:OAO:100216.100262} is known to have a best competitive ratio of $0.5$ for deterministic algorithms and $1-1/e$ for randomized algorithms. For the online weighted bipartite matching problem that we consider, the worst-case competitive ratio cannot be bounded below by any constant. Many subsequent works have tried to improve performance ratios under relaxed definitions of competitiveness.

Specifically, three types of assumptions are commonly used. The first type of assumption is that each demand node is independently and identically (i.i.d.) picked from a \emph{known} set of nodes. Under this assumption, \citet{doi:10.1287/moor.2013.0621, doi:10.1287/moor.1120.0551, Bahmani:2010:IBO:1888935.1888956, Feldman:2009:OSM:1747597.1748029} propose online algorithms with competitive ratios higher than $1-1/e$ for the cardinality matching problem, in which the goal is to maximize the total number of matched pairs. \citet{Haeupler:2011:OSW:2188743.2188759} study online algorithms with competitive ratios higher than $1-1/e$ for the weighted bipartite matching problem. Our definition of competitive ratio is the same as theirs. Our model is also similar, but we allow a more general arrival process of demand nodes in which the distribution of nodes can change over time.  Previous analyses depend crucially on the fact that demand nodes are i.i.d. in order to simplify the expression for the probability that any demand node is matched to any resource node.  The expression becomes much more complex, and the arguments break down in the case that demand arrivals are no longer i.i.d.


The second type of assumption is that the sequence of demand nodes is a random permutation of an \emph{unknown} set of nodes.  This random permutation assumption has been used in the secretary problem \citep{Kleinberg:2005:MSA:1070432.1070519,Babaioff:2008:OAG:1399589.1399596}, adword problem \citep{Goel:2008:OBM:1347082.1347189} and the bipartite matching problem \citep{Mahdian:2011:OBM:1993636.1993716,Karande:2011:OBM:1993636.1993715}. \citet{Kesselheim2013} study the weighted bipartite matching problem with extension to combinatorial auctions. Our work is different from all of these in that the non-stationarity of arrivals in our model cannot be captured by the random permutation assumption.

The third type of assumption made is that each demand node requests a very small amount of resource. The combination of this assumption and the random-permutation assumption often leads to polynomial-time approximation schemes (PTAS) for problems such as adword \citep{Devanur09theadwords}, stochastic packing \citep{Feldman:2010:OSP:1888935.1888957}, online linear programming \citep{doi:10.1287/opre.2014.1289}, and packing problems \citep{doi:10.1287/moor.2013.0612}. Typically, the PTAS proposed in these works use dual prices to make allocation decisions. Under this third assumption, \citet{Devanur:2011:NOO:1993574.1993581} study a resource allocation problem in which the distribution of nodes is allowed to change over time, but still needs to follow a requirement that the distribution at any moment induce a small enough offline objective value. They then study the asymptotic performance of their algorithm. In our model, the amount capacity requested by each customer is not necessarily small relative to the total amount of capacity available.  Therefore, the analysis in these previous works does not apply to our problem.


In our model, the arrival rates, or the distribution of demand nodes, are allowed to change over time. This non-stationarity poses new challenges, because it cannot be analyzed with existing methods.  At the same time, it is an essential feature in our model because it allows us to capture the perishability of service capacity in the applications that we consider.  When a resource perishes within the horizon, the demand for that resource drops to $0$.  Such a demand process must be time-varying. This important feature has received only limited attention so far. \citet{doi:10.1287/moor.1120.0548} consider an allocation model with a very general arrival process, but their allocation policy has performance guarantee only when the arrival rates are uniform.  In this paper, we allow arrival processes to be non-homogeneous Poisson processes with arbitrary rates.

Our algorithms solves a linear program and uses its optimal solution to make matching decisions. The idea of using optimal solutions to a linear program is natural and has been used by several previous works mentioned above. For example, \citet{Feldman:2009:OSM:1747597.1748029}, \citet{doi:10.1287/moor.1120.0551}, \citet{Haeupler:2011:OSW:2188743.2188759}, and \citet{Kesselheim2013} have used similar algorithms to obtain constant competitive ratios, albeit for different demand models.

The paper of \citet{alaei2012online} solves an online matching problem with non-stationary arrivals in a discrete-time setting.  They propose an algorithm similar to our Separation Algorithm, which is a primitive algorithm that we analyze initially and later improve upon.  They prove that this algorithm achieves a competitive ratio of at least $1-\frac{1}{\sqrt{k-3}}$ and at most approximately $1-\frac{1}{\sqrt{2\pi k}}$, where $k$ is the minimum capacity of a resource. Compared to \citet{alaei2012online}, we prove a stronger lower bound of $\ratiok$ on the competitive ratio for our Separation Algorithm, using a few of the same ideas but largely different techniques, as we will elaborate on in Section \ref{sec:newbound}.  Thus, our lower bound is more similar to their upper bound.  We also point out that the Separation Algorithm is not practical because it might route customers to resources that are already exhausted, while there are still other available resources.  More importantly, because of randomization, it might reject a high-priority customer, but accept a low-priority customer at nearly the same time.  For this reason, we propose a new ``bid-pricing'' algorithm, based on the Separation Algorithm, that avoids all of the above problems.  We prove that the improved algorithm has the same theoretical performance guarantee, and has much better computational performance as tested on real data.  


\subsection{Revenue Management}

Our work is also related to the revenue management literature. We refer to \citet{TalluriV2004} for a comprehensive review of this literature. Traditional works in this area assume that demands for products are exogenous and independent of the availability of other products \citep{Lautenbacher:1999:UMD:767713.768452,doi:10.1287/trsc.27.3.252,Littlewood}. The decision is whether to admit or reject a customer upon her arrival. Our model reduces to this admission control problem in the special case that the resources are identical and are available at the same time.

When customers are open to purchase one among a set of different resources, our model controls which resource to assign to each customer. Thus, our model captures the problem of managing opaque products. Sellers of an opaque product conceal part of the products' information from customers. Sellers have the ability to select which specific product to offer after the purchase of opaque product. This enables the seller to more flexibly manage their inventory. Opaque products are often sold at a discount compared to specific products, making them attractive to wider segments of the market. These products are common in internet advertising, tour operations, property management \citep{gallego2004managing} and e-retailing. Customers purchase an opaque product if the declared characteristics fit their preferences. The buyer agrees to accept any specific product that meets the opaque description.  In our model, a specific product corresponds to a node on the right side of a bipartite graph. A unit of demand for an opaque product corresponds to a node on the left that connects to all of the specific products contained in the opaque product. The weight of an edge corresponds to the revenue earned by selling the opaque product. 

Previous works related to opaque products include \citet{doi:10.1287/msom.1040.0054}, \citet{doi:10.1287/mksc.1070.0318}, 
\citet{Petrick20102027}, \citet{Chen2010897}, \citet{Lee20121641}, \citet{Gonsch2014280} and \citet{doi:10.1287/mnsc.2014.1948}. Due to the problem of large state space, most analyses focus on models with very few product types. For systems with many product types, some pricing and allocation heuristics are known.  There is numerical evidence that much of the benefit of opaque products can be obtained by having two or three alternatives \citep{elmachtoubW2013retailing}. However, when a retailer has a large number of alternative products, it is unclear how to design such an opaque product.  Our work focuses on online allocation policies with constant performance guarantees for the management of an opaque product with an arbitrary number of alternatives.

Our model assumes independent demands, i.e., the demand for each product is exogenous and independent of the availability of other products. Many recent works in revenue management consider endogenous demands, which means that customers who find their most preferred product unavailable might turn to other products.  
  Examples of works on dependent demands include \citet{gallego2004managing}, \citet{doi:10.1287/opre.1050.0194}, \citet{LiuV2008} and \citet{doi:10.1287/opre.2014.1328}. One of the main characteristics of these models is that customer preferences cannot be observed until purchase decisions are made. In such situation, sellers only have a distributional information of customer preferences. This phenomenon does not apply to admission scheduling systems. In these systems, customer preference can be revealed before a unit of a resource is assigned. In MyChart, for example, the system is able to customize the appointment to offer to each patient after knowing the patient's profile and availability. We assume that each customer's preference is observed before a resource is assigned. Knowledge of preferences gives service providers the ability to improve the efficiency of the resource allocation process by tailoring the service offered to each customer.


Our work is related to the still limited literature on designing policies for revenue management that are robust to the distribution of arrivals. \citet{Ball2009} analyze online algorithms for the single-leg revenue management problem. Their performance metric is the traditional competitive ratio that compares online algorithms with an optimal offline algorithm under the worst-case instance of demand arrivals. They prove that the competitive ratio cannot be bounded below by any constant when there are arbitrarily many customer types. In our work, we relax the definition of competitive ratio, and show that our algorithms achieve a constant competitive ratio (under our definition) for any number of customer types and for a more general multi-resource model.
\citet{CongApproximationRM} study approximation algorithms for an admission control problem for a single resource when customer arrival processes can be correlated over time.  They use as the performance metric the ratio between the expected cost of their algorithm and that of an optimal online algorithm. Our performance metric is stronger than theirs as we compare our algorithms against an optimal offline algorithm, instead of the optimal online policy. \citet{CongApproximationRM} prove a constant approximation ratio for the case of two customer types, and also for the case of multiple customer types with specific restrictions. They allow only one type of resource to be allocated. In our model, we assume arrivals are independent over time, but we allow for multiple customer types and multiple resources without additional assumptions. 

\section{Model}

Throughout this paper, we let $\bN$ denote the set of positive integers. For any $n \in \bN$, let $[n]$ denote the set $\{1,2,...,n\}$.

There are $n \in \bN$ resources and $m \in \bN$ customer types. Customers of each type $i \in [m]$ randomly arrive over a continuous horizon $[0,1]$ according to a known non-homogeneous Poisson process with rate $\lambda_i(t)$, for $t \in [0,1]$. Let $\Lambda_i \equiv \int_0^1 \lambda_i(t) \,dt$ be the expected total number of arrivals of type-$i$ customers.  Each resource $j \in [n]$ has a capacity of $C_j \in \bN$ units.

When a customer arrives, one unit of capacity of an available resource must be assigned to the customer, or the customer must be rejected. 
 A customer of type $i \in [m]$ earns a reward $r_{ij}$ if assigned to resource $j \in [n]$. The objective is to allocate the resources to the customers to maximize the expected total reward from all of the allocated resources.

	This model captures the expiration of resources in the following sense. Suppose we assign an expiration time $t_j \in [0,1]$ to each resource $j \in [n]$. Then, for any customer type $i \in [m]$ such that $\lambda_i(t) > 0$ for some  $t > t_j$, we require $r_{ij} = 0$. In this way, the reward from assigning resource $j$ to any customer who arrives after the expiration time $t_j$ is $0$.

\subsection{Definition of Competitive Ratios}

Let $\delta_i$ be the actual total number of arrivals of type $i$ customers. We must have $\mathbf{E}[\delta_i] = \Lambda_i$, for all $i \in [m]$.   An offline algorithm knows $\delta = (\delta_1,\delta_2,...,\delta_m)$ at the beginning of the horizon. Let OPT$(\delta)$ be the optimal offline reward given the number of arrivals $\delta$. Note that an optimal offline algorithm does not need to know the time of each arrival, as the algorithm essentially solves a maximum weighted matching problem, between the customers and resources. An online algorithm, however, does not know the entire sample path of future arrivals, but only knows the arrival rates $\lambda_i(t)$, $i \in [m]$. In this paper, we define the \emph{competitive ratio} as the ratio between the expected reward of an online algorithm and the expected reward of an optimal offline algorithm.

\begin{definition}
An online algorithm is \emph{$c$-competitive} if its total reward ALG satisfies
\[ \mathbf{E}[\text{ALG}] \geq c \,\mathbf{E}[\text{OPT}(\delta)],\]
where the expectation is taken over the sample path of customer arrivals (the random vector $\delta$ is determined by the sample path of arrivals).
\end{definition}

\subsection{Offline Algorithm and Its Upper Bound}


Before introducing our online algorithms, we first characterize an optimal offline algorithm and an upper bound on the optimal offline reward.

In the offline case, the total number of arrivals $\delta_i$ of each customer type $i$ is known, and the exact arrival time is irrelevant. Given the $\delta_i$'s, the maximum offline reward OPT$(\delta)$ can be computed by solving a maximum weighted matching problem, which can be formulated as the following LP:
\begin{align}
\begin{split}
\text{OPT}(\delta) = \max & \,\, \sum_{i \in [m]} \sum_{j \in [n]} x_{ij} r_{ij}\\
\mbox{s.t. }& \sum_{j=1}^n x_{ij} \leq \delta_i, \,\, \mbox{ for } i\in[m]\\
& \sum_{i=1}^m x_{ij} \leq C_j, \,\, \mbox{ for } j\in[n]\\
& x_{ij} \geq 0, \,\, \mbox{ for } i\in[m];\ j\in[n].
\end{split}\label{eq:LP1}
\end{align}
where the decision $x_{ij}$ is the number of type-$i$ customers who are assigned to resource $j$. 
Let $\bar x(\delta)$ be an optimal solution to this LP. Then OPT$(\delta) = \sum_{i=1}^m \sum_{j=1}^n r_{ij} \bar x_{ij}(\delta)$.

We are interested in finding an upper bound on the expected optimal offline reward $\mathbf{E}[\text{OPT}(\delta)]$. We next show that LP (\ref{eq:LP2}), which uses $\mathbf{E}[\delta]$ instead of $\delta$ as the total demand, gives such an upper bound:
\begin{align}
\begin{split}
\max & \,\, \sum_{i \in [m]} \sum_{j \in [n]} x_{ij} r_{ij}\\
\mbox{s.t. }& \sum_{j=1}^n x_{ij} \leq \Lambda_i, \,\, \mbox{ for } i\in[m]\\
& \sum_{i=1}^m x_{ij} \leq C_j, \,\, \mbox{ for } j\in[n]\\
& x_{ij} \geq 0.
\end{split}\label{eq:LP2}
\end{align}

\begin{theorem}
\label{thm:upperbound}
The optimal objective value of (\ref{eq:LP2}) is an upper bound on $\mathbf{E}[\text{OPT}(\delta)]$.
\end{theorem}

\proof{Proof.}
Since $\sum_{j=1}^n \bar x_{ij}(\delta) \leq \delta_i$ and $\sum_{i=1}^m \bar x_{ij}(\delta) \leq C_j$, we must have $\sum_{j=1}^n \mathbf{E}[\bar x_{ij}(\delta)] \leq \mathbf{E}[\delta_i] = \Lambda_i$ and $\sum_{i=1}^m \mathbf{E}[\bar x_{ij}(\delta)] \leq C_j$. Thus, $\mathbf{E}[\bar x_{ij}(\delta)]$ is a feasible solution to the LP (\ref{eq:LP2}). It follows that the optimal objective value of (\ref{eq:LP2}) is an upper bound on \[\sum_{i=1}^m \sum_{j=1}^n r_{ij} \mathbf{E}[\bar x_{ij}(\delta)] = \mathbf{E}[\text{OPT}(\delta)].\]  
Similar techniques have been used in revenue management to prove similar results \citep{GallegoV1997}.
\halmos
\endproof


\section{Separation Algorithm}

In this section, we propose the {\it Separation Algorithm}. The algorithm works by solving the LP (\ref{eq:LP2}) once, routing the customers to the resources according to an optimal solution to the LP (\ref{eq:LP2}).  Then, for each resource separately, the algorithm optimally controls the admission of customers who have been routed to that resource.  Using the LP information with respect to the expected number of arrivals (or sometimes, an estimate of the expected number of arrivals) is natural and has been used in several previous results (for example,  \citet{Feldman:2009:OSM:1747597.1748029}, \citet{doi:10.1287/moor.1120.0551}, \citet{Haeupler:2011:OSW:2188743.2188759}, and \citet{Kesselheim2013}).


Let $x^*$ be an optimal solution to the linear program (\ref{eq:LP2}). Whenever a customer of type $i$ arrives, the Separation Algorithm randomly and independently picks a candidate resource $j \in [n]$ with probability $x^*_{ij} / \Lambda_i$, regardless of the availability of resources. We say that this customer is \emph{routed} to resource $j$.  
According to the Poisson thinning property, the arrival process of type-$i$ customers who will be routed to resource $j$ is a non-homogeneous Poisson process with rate
\begin{equation} \lambda_{ij}(t) \equiv \lambda_i(t) x_{ij}^*/\Lambda_i,\,\,\, \mbox{ for } 0\leq t \leq 1. \label{eq:arrivalrate}\end{equation} 

Viewing the random routing process as exogenous, each resource $j$ receives an independent arrival process with split rate $\lambda_{ij}(t)$ from each customer type $i$. Then for each resource $j$, the Separation Algorithm optimally controls the admission of customers who will be routed to resource $j$. That is, when a type-$i$ customer is routed to resource $j$ at time $t$, the algorithm compares $r_{ij}$ with the marginal cost of taking one unit away from resource $j$, where the marginal cost is computed based on the future customers who will be routed to resource $j$. The customer is accepted and offered resource $j$ if $r_{ij}$ is higher than or equal to the marginal cost. The customer is rejected if $r_{ij}$ is smaller than the marginal cost or if resource $j$ has no remaining capacity.

For each resource $j \in [n]$, let $c_j(t) \in \{0,1,...,C_j\}$ denote the amount of resource $j$ that remains at time $t$. Given the exogenous random routing process, we define $f_j(t,c_j(t))$ as the optimal future reward of the admission control problem for resource $j$. $f_j(t,c_j(t))$ is governed by the Hamilton-Jacobi-Bellman equation
\begin{equation}  \frac{\partial {f}_j(t,c)}{\partial t} = -\sum_{i \in [m]}  \lambda_{ij}(t) (r_{ij} -f_j(t,c) + f_j(t,c-1))^+, \ \ \forall c=1,2,...,C_j, \ \ t \in [0,1]. \label{eq:dpcapacity}\end{equation}
The boundary conditions are $f_j(1,c) = 0$ for all $c=0,1,...,C_j$, and $f_j(t,0) = 0$ for all $t \in [0,1]$. According to properties of the HJB equation, $f_j(t,c)$ is decreasing in $t$, which captures the fact that resources are expiring over time. We call $f_j(\cdot,\cdot)$ the \emph{reward function} for resource $j$. In practice, the continuous-time dynamic programming \eqref{eq:dpcapacity} can often be solved by discretizing the horizon (for example, see \citet{ContinuousRM}).

Below are the detailed steps of the Separation Algorithm:
\begin{enumerate}
\item Solve LP (\ref{eq:LP2}). Let $x^*$ be any optimal solution.
\item For each resource $j \in [n]$, compute the reward function $f_j(t,c)$ according to (\ref{eq:dpcapacity}).
\item Upon an arrival of a type-$i$ customer at time $t$, randomly pick a number $j \in [n]$ with probability $x_{ij}^* / \Lambda_i$. Assign resource $j$ to the customer if resource $j$ has positive remaining capacity and $r_{ij} \geq f_j(t,c_j(t)) - f_j(t,c_j(t)-1)$.
\end{enumerate}

The following proposition states that the total expected reward of the Separation Algorithm is given by the reward functions. We omit the proof as this result directly follows from properties of the HJB equation.
\begin{proposition}\label{prop:int}
Conditioned on the state $c_j(t)$, the Separation Algorithm earns reward $f_j(t,c_j(t))$ from resource $j$ in time $[t,1]$ in expectation. In particular, the expected total reward of the Separation Algorithm is $\sum_{j \in [n]} f_j(0,C_j)$. 
\end{proposition}


\section{Proof of Competitive Ratio} \label{sec:newbound}

In this section, we show that if $k$ is the minimum capacity of any resource, then the competitive ratio of the Separation Algorithm is $\ratiok$. This result is stated in Theorem \ref{thm:betaBound}. 

To prove the competitive ratio, we fix a resource $j \in [n]$ and  focus on the ratio
\begin{equation}\label{eq:ratioInitial}\frac{f_j(0,C_j)}{\sum_{i \in [m]} x^*_{ij} r_{ij}},
\end{equation}
where $f_j(0,C_j)$ is the expected reward that the Separation Algorithm earns from resource $j$, and $\sum_{i \in [m]} x^*_{ij} r_{ij}$ is an upper bound on the optimal expected offline reward from resource $j$ (see LP \eqref{eq:LP1}).



We want to lower-bound \eqref{eq:ratioInitial} by examining all possible inputs $r$ and $\lambda(\cdot)$. As we will prove a lower bound that increases in the capacity value, the worst case for our analysis is $C_j = k$, where $k$ is the minimum capacity of any resource. Thus, for the rest of the section, we suppose $C_j = k$ and analyze the ratio $\frac{f_j(0,k)}{\sum_{i \in [m]} x^*_{ij} r_{ij}}$.

We apply Jensen's inequality to the HJB equation \eqref{eq:dpcapacity} to obtain
\begin{align*}
\frac{\partial {f}_j(t,c)}{\partial t} & = -\sum_{i \in [m]}  \lambda_{ij}(t)  (r_{ij} -f_j(t,c) + f_j(t,c-1))^+\\
& \leq - \sum_{i \in [m]}  \lambda_{ij}(t)    \left( \frac{\sum_{i \in [m]}  \lambda_{ij}(t)r_{ij} }{\sum_{i \in [m]}  \lambda_{ij}(t)}  -f_j(t,c) + f_j(t,c-1)\right)^+.
\end{align*}
Thus, the performance of the Separation Algorithm can be \emph{lowered} by replacing the problem instance with one in which there is only one type of customer arrival with arrival rate $\lambda(t)=\sum_{i \in [m]}  \lambda_{ij}(t)$ and reward rate $r(t)= \frac{\sum_{i \in [m]}  r_{ij} \lambda_{ij}(t)}{\sum_{i \in [m]}  \lambda_{ij}(t)}$,
so that the worst-case instance has one customer type, and time-dependent reward value $r(\cdot)$.  This observation has also been made by \citet{alaei2012online}. 

Furthermore, by \eqref{eq:arrivalrate} and definition of $\Lambda_i$, we can obtain 
\begin{align*}
\sum_{i \in [m]} x^*_{ij} r_{ij} &=  \sum_{i \in [m]} r_{ij}\int_0^1 \lambda_{ij} (t)dt \frac{\int_0^1 \lambda_{i} (t)dt}{\Lambda_i} \\
& = \int_0^1 \sum_{i \in [m]} r_{ij} \lambda_{ij} (t)dt\\
&= \int_0^1 r(t) \lambda(t) dt.
\end{align*}

Thus, to characterize the worst-case performance ratio for the fixed resource $j$, we only need to lower-bound
\begin{equation}\label{eq:largecapacityratio} 
\frac{f_j(0,k)}{\sum_{i \in [m]} x^*_{ij} r_{ij}} \geq \frac{u_0(0)}{\int_0^1 r(t)\lambda(t) dt},
\end{equation}
where $u_l(t)$ is the new reward function defined based on $\lambda(t)$ and $r(t)$:
\begin{equation}\label{eq:hjbu}
 \frac{d u_l(t)}{d t} = -\lambda(t)(r(t)-u_l(t)+u_{l+1}(t))^+,\ \ \ \ \forall l =0,1,...,k-1,\ \ t \in [0,1]
 \end{equation}
with boundary conditions $u_l(1)=0$ for all $l=0,1,...,k$ and $u_k(t) = 0$ for all $t \in [0,1]$. 

Note that the HJB equation \eqref{eq:hjbu} is different from \eqref{eq:dpcapacity}, as \eqref{eq:hjbu} is defined by a different arrival rate and reward rate. Moreover, we use $l$ to denote the \emph{consumed inventory} in \eqref{eq:hjbu}, whereas in \eqref{eq:dpcapacity}, we used $c$ to denote the \emph{remaining inventory}. This change is convenient for our analysis.   

In order to lower-bound \eqref{eq:largecapacityratio}, we need to examine all possible reward rates $r(\cdot)$ and arrival rates $\lambda(\cdot)$ such that the constraints of \eqref{eq:LP2} are satisfied. The first constraint of \eqref{eq:LP2} is satisfied by definition of $\lambda(t)$ and $\lambda_{ij}(t)$. The second constraint of (\ref{eq:LP2}) requires 
\begin{equation}\label{eq:capacityinequality} \int_0^1 \lambda(t) dt = \int_0^1 \sum_{i \in [m]}\lambda_{ij}(t)dt  =\sum_{i \in [m]} x^*_{ij}  \leq k.
\end{equation}


\subsection{Homogenizing time}
Without loss of generality, we can change the horizon length, the arrival rate $\lambda(\cdot)$ and the reward rate $r(\cdot)$ as follows, while keeping the ratio (\ref{eq:largecapacityratio}) unchanged:
 \begin{enumerate}
\item If the inequality \eqref{eq:capacityinequality} is not tight, i.e., $\int_0^1 \lambda(t) dt < k$, we can extend the horizon to  length $T>1$ by adding more arrivals with reward $0$. Thus, we can equivalently assume $\int_0^T \lambda(t) dt = k$.
\item Let $T$ be the length of the (possibly extended) horizon such that $\int_0^T \lambda(t) dt = k$. Define a virtual time
\[ \bar t(t) \equiv  \int_0^t \lambda(s) ds, \ \ \forall t \in [0,T].\]
We must have $\bar t(t) \in [0,k]$ for all $t \in [0,T]$. Using this new time variable $\bar t(\cdot)$, we can define new reward functions as
\[ \bar u_l(s) = u_l(\bar t^{-1} (s)),\]
where we interpret $\bar t^{-1}(s)$ as the first time $t$ that satisfies $\bar t(t) = s$. Similarly, we can define $\bar r(s) = r(\bar t^{-1}(s))$. Then we can equivalently transform the HJB equation for $u_l(t)$ as follows
\[ \frac{d u_l(t)}{d t} =\frac{d \bar u_l(\bar t(t))}{d \bar t(t)}\frac{d \bar t(t)}{dt} = \frac{d \bar u_l(\bar t)}{d \bar t} \lambda(t)\]
\[ \Longrightarrow   \frac{d \bar u_l(\bar t)}{d \bar t} \lambda(t) =- \lambda(t)(r(t) + u_{l+1}(t) - u_l(t))^+ = - \lambda(t)(\bar r(\bar t) + \bar u_{l+1}(\bar t) - \bar u_l(\bar t))^+  \]
\[ \Longrightarrow \frac{d \bar u_l(\bar t)}{d \bar t} =  -(\bar r(\bar t) + \bar u_{l+1}(\bar t) - \bar u_l(\bar t))^+, \,\,\,\forall \bar t \in [0,k].\]
This equation can be viewed as another HJB equation with arrival rate $1$ and reward rate $\bar r(\cdot)$, with boundary conditions $u_{k}(t)=0$ for $t\in[0,k]$ and $u_l(k)=0$ for $l=0,1,...,k$. Furthermore, the upper bound on the expected offline reward can be transformed as
\[ \int_0^T r(t) \lambda(t) dt = \int_0^{k} \bar r(\bar t) d\bar t.\]
In summary, we can equivalently transform the problem into one whose arrival rate is uniformly $1$ and whose time horizon is $[0,k]$.
\end{enumerate}

\subsection{Bound-revealing optimization problem}

After applying the above transformations, we can write an optimization problem that reveals the competitive ratio as follows
\begin{align}\label{boundRevealingProblem}
\min_{r(t), u_i(t), i=0,1,...,k-1; t \in [0,k]} & u_0(0)\\
\nonumber\text{s.t. } & \frac{d u_i(t)}{d t} = - (r(t) + u_{i+1}(t) - u_i(t))^+, \,\,\forall i = 0,1,...,k-1; t \in [0,k]\\
\nonumber& \int_0^k r(t) dt = 1\\
\nonumber& u_i(t) \geq 0, \,\, \forall i = 0,1,...,k-1; t\in[0,k]\\
\nonumber&  r(t) \geq 0.
\end{align}

Here the second constraint  $\int_0^k r(t) dt = 1$ normalizes the upper bound on the expected offline reward. By using $g_i(t) = - d u_i(t) / d t$ and replacing $(\cdot)^+$ with linear constraints, we can write the above problem equivalently as (note that $g_k(t) = 0, \forall t \in [0,k]$)
\begin{align*}
\min_{r(t),g_i(t), i = 0,1,...,k-1; t \in [0,k]} & \int_0^k g_0(s) ds \\
\text{s.t. } &  g_i(t) \geq r(t) + \int_t^k g_{i+1}(s) ds -  \int_t^k g_i(s) ds, \,\,\forall i = 0,1,...,k-1; \forall t \in [0,k]\\
& \int_0^k r(t) dt = 1\\
& g_i(t) \geq 0, \,\,\forall i = 0,1,...,k-1 ; t \in[0,k]\\
&  r(t) \geq 0.
\end{align*}

Let $\alpha_i(t)$ be a dual variable for the first constraint, for all $i=0,1,...,k-1$ and $t \in [0,k]$. Let $\beta$ be a dual variable for the second constraint. The dual problem is

\begin{align}
\begin{split}\label{eq:dual}
\max_{\alpha_i(t), \beta} & \,\,\, \beta\\
\text{s.t. } & \alpha_0(t) + \int_0^t \alpha_0(s) ds \leq 1, \,\,\forall t \in [0,k]\\
& \alpha_i(t) + \int_0^t \alpha_i(s) ds \leq \int_0^t \alpha_{i-1}(s)ds , \,\,\forall i=1,2,...,k-1; \forall t \in [0,k] \\
& \beta \leq \sum_{i=0}^{k-1} \alpha_i(t) , \,\,\forall t \in [0,k]\\
& \alpha_i(t) \geq 0, \,\,\forall i=0,2,...,k-1, \forall t \in [0,k].
\end{split}
\end{align}

This dual problem tries to maximize the minimum value of $\sum_{i=0}^{k-1} \alpha_i(t)$ with respect to $t$. The optimal $\beta$ is a lower bound on the competitive ratio that we seek to characterize.

\subsection{A dual-feasible solution for the bound-revealing problem}

We first show that a feasible solution to the dual problem (\ref{eq:dual}) can be constructed based on a modification of a Poisson process.  As we shall explain shortly, this is a Poisson process to which we apply a control, using a sequence of \emph{bounding barriers}.  We will use the solution obtained via this derived process to obtain a lower bound on the optimal value of the bound-revealing optimization problem \eqref{boundRevealingProblem}. We will refer to the process as a \emph{bounded} Poisson process.  \citet{alaei2012online} also prove their bound by working with a dual-feasible solution.  However, we construct our dual-feasible solution differently using a novel method.  Because our bound has to be tighter, our analysis of this solution is also much more involved.

Let $t_0, t_1,t_2,...,t_k$ be a sequence of time points such that $0=t_0 < t_1 < \cdots < t_{k-1} < t_k = k$.


\begin{figure}[h!]
   \centering
   \caption{Illustration of the bounded Poisson process. Dashed red line is the barrier. Solid blue line is the $R(t)$ process. The barrier is active when the two lines overlap.}
   \label{fig:boundedprocess}
     \includegraphics[width=0.5\textwidth]{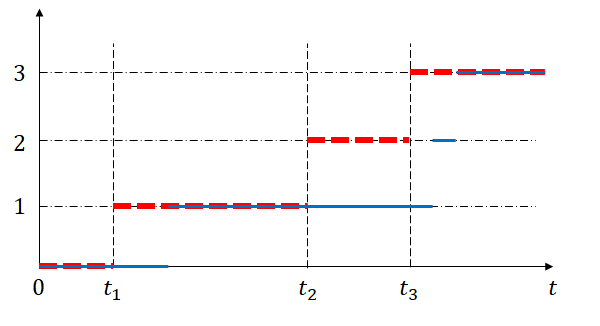}
\end{figure}

Let $\{N(t)\}_{t\geq 0}$ be a (counting) Poisson process with rate $1$. We apply an upper barrier to $N(t)$ to obtain a new bounded process $\{R(t)\}_{t\geq 0}$. Figure \ref{fig:boundedprocess} illustrates this $R(t)$ process. Starting with an initial value $0$ at time $t_0=0$, the barrier increases by $1$ at times $t_1,t_2,...,t_{k-1}$. At these time points, the new bounded process has values
\[ R(t_i) = \min(i-1, R(t_{i-1}) + N(t_i) - N(t_{i-1})), \,\forall i = 1,2,...,k-1,\]
with $R(t_0) = R(0) = 0$.  And for $t \in [t_i, t_{i+1}]$, we have
\[ R(t) = \min(i, R(t_i) + N(t) - N(t_i)), \,\forall i = 1,2,...,k-1.\]

Eventually,
\[ R(t_k) = R(k) = \min(k-1, R(t_{k-1}) + N(k) - N(t_{k-1})).\]


\begin{theorem}\label{thm:nboundthm1}
There exists a feasible dual solution $\beta^*$, $\alpha^*_i(t)$ for $t\in[0,k]$, $i=0,1,2,...,k-1$, such that 
\[ \alpha_i^*(t) = \bP(R(t) = i) \beta^*, \,\,\,\forall t \in[0,k], i=0,1,...,k-1,\]
\begin{equation}\label{eq:nboundthm1eq1} k (1- \beta^*) = \beta^*\left[k - \sum_{i=0}^{k-1} i \bP(R(k) = i)\right] \end{equation}
for the bounded Poisson process $R(t)$ as constructed above.
\end{theorem}

Given that $\beta^*$ and $\alpha^*$ are dual-feasible, we will next attempt to lower-bound objective $\beta^*$ by analyzing the process $R(\cdot)$.  

First we show that the times at which the barriers are applied are bounded above by $1, 2,\ldots,k-1$.
\begin{theorem}\label{thm:nboundlmtp}
The time points $t_1, t_2,...,t_{k-1}$ constructed in the proof of Theorem \ref{thm:nboundthm1} satisfy $t_i \leq i$, for $i=1,2,...,k-1$.
\end{theorem}

Before proving Theorem \ref{thm:nboundlmtp}, we first prove Lemmas \ref{lm:nboundlm4} to \ref{lm:nboundlmkey}, which further characterize the behavior of the process $R(t)$.  These lemmas collectively show that when the barriers are applied at regular points starting at some time of the horizon, i.e., $t_i = i$ for all $i\geq l$ for some integer $l$, the time spent at the barriers must be monotone decreasing in the index $i$ for all $i \geq l$.

For ease of notation, let
\[ I_i \equiv \int_{t_i}^{t_{i+1}} \mathbf{1}(R(s) = i) ds\]
be the total time that the bounded process $R(t)$ stays at the barrier $i$ during the interval $[t_i, t_{i+1}]$, for $i=0,1,...,k-1$. Note that $\mathbf{E}[I_i] = \int_{t_i}^{t_{i+1}} \bP(R(s) = i) ds$.  Let $P_i(\lambda)$ be the probability that a Poisson random variable with mean $\lambda$ is equal to $i$. Let $P_{\geq i}(\lambda)$ and $P_{\leq i}(\lambda)$ denote $\sum_{j=i}^\infty P_j(\lambda)$ and $\sum_{j=0}^i P_j(\lambda)$, respectively.

First, assuming that the barriers are applied at regular points $0,1,\ldots,k-1$, we can quantify the difference in the expected time spent at each barrier, given different starting points for the process $R(\cdot)$.
\begin{lemma}\label{lm:nboundlm4}
Given any $l \in \{1,2,...,k-1\}$, if $t_l = l$ and $t_{l+1} = l+1$, we must have 
\[ \mathbf{E}[I_l | R(l) = l-j] - \mathbf{E}[I_l | R(l) = l-j-1] = P_{\geq j+1}(1) \]
for all $j=0,1,...,l-1$.
\end{lemma}

Next, assuming that the barriers are applied at regular points $0,1,\ldots,k-1$, we can characterize the differences in the expected time spent at each barrier for successive pairs of starting points.
\begin{lemma}\label{lm:nboundlmkeybound}
Given any $l \in \{2,3,...,k-1\}$, if $t_i = i$ for all $i=l,l+1,...,k-1$, we must have
\[ \mathbf{E}[I_i|R(l) = l] - \mathbf{E}[I_i | R(l) = l-1] \geq e^{-1}(\mathbf{E}[I_i|R(l) = l-1] - \mathbf{E}[I_i | R(l) = l-2])\]
for all $i = l,l+1,...,k-1$.
\end{lemma}

Using the previous result, we relax the assumption that all the barriers are applied at regular points $0,1,\ldots,k-1$.  We assume now that the barriers are applied at regular times beyond a point.  Under this condition, we show that the differences in the expected time spent at successive barriers are increasing with the starting point of the process.
\begin{lemma}\label{lm:nboundlm5}
Given any $l \in \{1,2,...,k-2\}$, if $t_i \leq i$ for $i=1,2,...,l$, and $t_i = i$ for $i=l+1,l+2,...,k-1$, we must have 
\[ \mathbf{E}[I_i | R(l) = l] - \mathbf{E}[I_{i+1}|R(l) = l] \geq \mathbf{E}[I_i | R(l) = l-1] - \mathbf{E}[I_{i+1} | R(l) = l-1]\]
for all $i = l,l+1,...,k-2$.
\end{lemma}

Next, assuming that the barriers are applied at regular points $0,1,\ldots,k-1$, we show that the expected time spent by the process at each barrier is decreasing with the index of the barrier.
\begin{lemma}\label{lm:nboundlmmonotone}
If $t_i = i$ for all $i=1,2,...,k-1$, we must have $\mathbf{E}[I_i] \geq \mathbf{E}[I_{i+1}]$ for all $i=1,2,...,k-2$.
\end{lemma}
\proof{Proof.}
It is obvious that for any $i \geq 1$,
\[ \mathbf{E}[I_i | R(1) = 1] \geq \mathbf{E}[I_i | R(1) = 0],\]
because when the starting position becomes lower, it is harder for the random process $R(t)$ to reach the barrier at any later time. Since $\mathbf{E}[I_i | R(1) = 0] = \mathbf{E}[I_i]$, and by symmetry, $\mathbf{E}[I_i | R(1) = 1] = \mathbf{E}[I_{i-1}]$, we have $\mathbf{E}[I_{i-1}] \geq \mathbf{E}[I_i]$ for all $i \geq 1$.
\halmos\endproof

Finally, we relax the requirement of Lemma \ref{lm:nboundlmmonotone}.  We require only that the barriers be applied at regular points only after some time.  We show that the expected time spent at the barriers are still decreasing.
\begin{lemma}\label{lm:nboundlmkey}
Given any $l \in \{1,2,...,k-2\}$, if $t_i \leq i$ for $i=1,2,...,l$, and $t_i = i$ for $i=l+1,l+2,...,k-1$, we must have 
\[\mathbf{E}[I_i] \geq \mathbf{E}[I_{i+1}]\]
for all $i=l,l+1,...,k-2$.
\end{lemma}

The idea of the proof of Theorem \ref{thm:nboundlmtp} is as follows.  We will start by setting the barriers at times $0,1,\ldots,k-1$.  We then successively reduce the values $t_i$, $i=0,1,\ldots$, until the expected time spent at each barrier is no more than $1/\beta-1$.  By the monotonicity shown in Lemma \ref{lm:nboundlmkey}, this procedure must stop with the expected time spent at each barrier bounded above by $1/\beta-1$.

If we change the value of $\beta$, the time points $t_1,t_2,...,t_{k-1}$ that result from the above procedure must change continuously in $\beta$. We simply choose $\beta$ such that, when the procedure ends, the expected time spent at the last barrier is $1/\beta-1$, which implies that the expected time spent at all barriers is exactly $1/\beta-1$.

\subsection{Computing the bound}
First, we prove an inequality, which will be useful in computing our bound.
\begin{lemma}\label{lm:nboundlm3}
 For any $x,y \in \mathbb{Z}$ and $\lambda \in [0,k]$ such that $x\geq y \geq k-1-\lambda$, we must have for any $l = 0,1,...,k-1$,
\[ \sum_{i=-l}^l P_{k-1+i-x}(\lambda) \leq \sum_{i=-l}^l P_{k-1+i-y}(\lambda).\]
\end{lemma}

Finally, we derive our lower bound on $\beta^*$.  The bound is simply a reduction of the equation
\[ k (1- \beta^*) = \beta^*\left[k - \sum_{i=0}^{k-1} i \bP(R(k) = i)\right], \]
which follows from Theorem \ref{thm:nboundthm1}.  $\beta^*$ is strictly greater than $0.5$ for $k\geq 2$.  For example, when $k = 2$, $\beta^*$ satisfies
\[ 3 \beta + \beta e^{1/\beta - 3} = 2,\]
from which we can obtain $\beta^* \approx 0.615$. 
\begin{proposition} \label{prop:nboundthmy}
\[ \beta^* \geq \frac{1}{1 + \frac{1}{k} [ \sum_{i=2k-1}^\infty i P_i(k) +2 \sum_{i=1}^{k-1} i P_{k+i-1}(k)]}. \]
\end{proposition}

\begin{theorem}\label{thm:betaBound} If $k$ is the minimum capacity of any resource, the competitive ratio for the Separation Algorithm is at least
\[ \beta^* \geq \frac{1}{1 + 2 \left[\frac{P_{\geq k}(k)}{k} + \frac{e^{-k} k^k}{k!}\right]} = 1 - \sqrt{\frac{2}{\pi}} \frac{1}{\sqrt{k}} + O(\frac{1}{k}).\]
\end{theorem}

\section{Marginal Allocation Algorithm}

In this section, we propose the Marginal Allocation Algorithm, which improves on the Separation Algorithm by converting it to a bid-price algorithm. 

The Separation Algorithm, when carried out in practice, has several problems.  First, it might route customers to unavailable resources when they can be better matched to other resources.  Second, because of the random routing, it might unfairly accept a lower-priority customer after rejecting a higher-priority customer.  In this section, we present the Marginal Allocation Algorithm\ which resolves these issues by converting the Separation Algorithm into a bid-price algorithm. We will prove that the Marginal Allocation Algorithm\ has theoretical performance no worse than that of the Separation Algorithm. 

The Marginal Allocation Algorithm\  uses the marginal reward $f_j(t,c_j(t)) - f_j(t,c_j(t)-1)$ as a bid price for resource $j$. When a customer of type $i$ arrives, the Marginal Allocation Algorithm\  rejects the customer if $r_{ij} < f_j(t,c_j(t)) - f_j(t,c_j(t)-1)$ for all available resources $j$; otherwise, it assigns this customer to resource
\[ \argmax_j\{ r_{ij} - f_j(t,c_j(t)) + f_j(t,c_j(t)-1) | c_j(t) > 0\} .  \]

To carry out this algorithm, we only need to compute the $n$ reward functions at the beginning of the horizon. Thus the space requirement is polynomial in $T$ and $n$. At any time $t$, we only need to know the $n$ reward functions $f_j(t,c_j(t))$, for $j = 1,2,...,n$, so as to make a decision.
 
The following theorem states that the Marginal Allocation Algorithm\  performs at least as well as the Separation Algorithm:
\begin{theorem}\label{thm:algdomination}
The expected total reward of the Marginal Allocation Algorithm is no less than that of the Separation Algorithm.
\end{theorem}
%
%

As a result, when $k$ is the minimum capacity, the competitive ratio of the Marginal Allocation Algorithm  is $\ratiok$. When $k$ tends to infinity, the competitive ratios tends to $1$, so the Marginal Allocation Algorithm is asymptotically optimal.

\subsection{Resource sharing}

In settings in which customers have similar preferences for all the resources, the Marginal Allocation Algorithm utilizes resources more effectively than the Separation Algorithm, because the latter restricts each customer to only one resource but the former can allocate any available resource. We focus on such settings in this section, and lower-bound the expected reward that the Marginal Allocation can earn more than the Separation Algorithm.

\begin{proposition}
The expected reward earned by the Marginal Allocation Algorithm can be as much as $1/(1-e^{-1})$ times that earned by the Separation Algorithm.
\end{proposition}
\proof{Proof.}
Suppose $r_{ij} = 1$ for all $j \in [n]$ and $i \in [m]$. Suppose $C_j = 1$ for all $j \in [n]$. Suppose $\sum_{i \in [m]} \Lambda_i = n$. In this way, the total expected number of arrivals is equal to the total capacity.

The optimal LP solution must satisfy $\sum_{i \in [m]} x_{ij}^* = 1$.

Since for any resource $j \in [n]$, the reward values $r_{ij}$ are the same for all customer types $i \in [m]$, the expected future reward earned from future customers must be no more than the reward values, i.e., $f_j(t,1) \leq r_{ij}= 1$ for all $t \in [0,1]$. 
\begin{align*}
f_j(0,1) & = \int_0^1 \sum_{i \in [m]} \lambda_{ij}(s) (r_{ij} - f_j(s,1) + f_j(s,0))^+ ds\\
& = \int_0^1 \sum_{i \in [m]} \lambda_{ij}(s) (1 - f_j(s,1))^+ ds\\
& = \int_0^1 \sum_{i \in [m]} \lambda_{ij}(s) (1 - f_j(s,1)) ds.
\end{align*}
\[ \Longrightarrow f_j(0,1) = 1 - e^{-\int_0^1 \sum_{i \in [m]} \lambda_{ij}(s) ds} = 1 - e^{-\int_0^1\sum_{i \in [m]} \lambda_i(s) x_{ij}^*/\Lambda_i ds} = 1-e^{-1}.\]
This result is easy to see, as the Separation Algorithm collects a reward the first time a customer is routed to resource $j$. $1-e^{-1}$ is exactly the probability that the number of customers routed to resource $j$ (Poisson random variable with mean $1$) is greater than or equal to $1$.

Thus, the total expected reward earned by the Separation Algorithm is
\begin{align*}
\sum_{j \in [n]} f_j(0,1) = n (1 - e^{-1}).
\end{align*}

In this setting, since the bid prices are no larger than the reward values, the Marginal Allocation Algorithm never rejects a customer unless all the resources are empty. Let $N$ be the total number of arrivals during the horizon. $N$ is a Poisson random variable with mean $\sum_{i \in [m]} \Lambda_i = n$.

According to Chebyshev's inequality,
\[ \bP( N < \bE[N] - (\bE[N])^{0.75}) \leq \frac{\Var(N)}{(\bE[N])^{1.5}} = n^{-0.5}. \]

The Marginal Allocation Algorithm earns at least
\begin{align*}
& [\bE[N] -  (\bE[N])^{0.75}]\bP(N \geq \bE[N] - (\bE[N])^{0.75})\\
 = & [n -  n^{0.75}][1 - \bP(N < \bE[N] - (\bE[N])^{0.75})]\\
 \geq & [n -  n^{0.75}][1 - n^{-0.5}]\\
 \geq & n - n^{0.75} - n^{0.5}.
\end{align*}

 The ratio between the total expected reward earned by the Marginal Allocation Algorithm and the Separation Algorithm is at least
\[ \frac{n - n^{0.75} - n^{0.5}}{n(1-e^{-1})} = \frac{1}{1-e^{-1}}(1 - n^{-0.25} - n^{-0.5}).\]
It approaches $1/(1-e^{-1})$ when $n$ becomes large.

\halmos
\endproof

\subsection{Overbooking}
\label{sec:overbooking}
No-shows is an issue that is common to all advance admission-scheduling systems.  When customers book in advance, events may transpire between the date of the booking and the planned date of service that cause customers to miss their appointments. Due to the frequent occurrence of no-shows, overbooking is commonly used in service industries. Suppose each customer has a no-show probability of $p_j$ when assigned to resource $j$, and incurs a cost of $D_j$ when being denied getting  resource $j$. Then we can model the overbooking strategy by expanding capacities at additional costs. Assume that the no-show events are exogenous to both online and offline algorithm. For resource $j$, the $k$th overbooked unit of capacity incurs an expected marginal cost of
\begin{equation}\label{eq:overbookingcost}
 o_j(k) = D_j \cdot (1- p_j) \cdot \left[ \sum_{l=0}^{k-1} \left(\!\!\begin{array}{cc}C_j+k-1\\l\end{array}\!\!\right) p_j^{l}(1-p_j)^{C_j + k-1-l} \right],
 \end{equation}
where the value in the brackets represents the probability that, among the $C_j + k-1$ customers who have already booked resource $j$, at most $k-1$ of them do not show up. The additional $1-p_j$ in the product represents the probability that the $k$th overbooked customer does show up. Note that the marginal cost $o_j(k)$ is independent of customer type, and is increasing in $k$.

Assuming that the reward $r_{ij}$ is earned whether a customer of type $i$ actually takes resource $j$, the marginal reward of allocating the $k$th overbooked unit of resource $j$ to a type $i$ customer is
\[ \tilde r_{i,j,k} = r_{ij} - o_j(k).\]

When using this reward value $\tilde r_{i,j,k}$, we are treating each overbooked unit of resource $j$ as a virtual slot to be allocated. Then, the theoretical bound of our algorithms still applies, with $\tilde r_{i,j,k}$ being the reward of expanded units.

Since $\tilde r_{i,j,k} \leq r_{i,j}$ and $\tilde r_{i,j,k}$ decreases in $k$, an optimal offline algorithm, when allocating resource $j$, will first fill in the $C_j$ units of regular capacity and then assign customers to those virtual slots with lower values of $k$. It will not use virtual slots with non-positive marginal reward. Then, when $b$ overbooked units of resource $j$ are used under the optimal offline algorithm by the end, the total cumulative cost
\begin{equation} \sum_{k=1}^b o_j(k)\label{eq:overbookingCost}\end{equation}
is just the actual expected overbooking cost for resource $j$.

\section{Numerical Studies}
\label{sec:numerical}

We compare our Marginal Allocation Algorithm against the outcome of the actual scheduling practices used in the Division of Clinical Genetics within the Department of Pediatrics at Columbia University Medical Center (CUMC). The third author oversees appointment scheduling practice at the medical center.  We estimate our model parameters, including patient preferences, arrival rates and hospital processing capacities, by using historical appointment-scheduling data from the outpatient clinics. 
We also test the performance of our algorithm against some simple heuristics. We find that our Marginal Allocation Algorithm performs the best among all heuristics considered, and is $21\%$ more efficient than current practice, according to our performance metric, which we will explain below.

Specifically, we used data from the Division of Clinical Genetics at CUMC.  Clinical Genetics is a field of medicine where adults are assessed for the risk of having offsprings with heritable conditions and children are assessed for genetic disorders.  Geneticists use physical exams, chromosome testing and DNA analysis to diagnose patients suspected of having genetic abnormalities. The data we used contain more than 9000 appointment entries recorded in the year 2013. Each entry in the data records information about one appointment.  The entry includes the date that the patient makes the appointment, the exact time of the appointment, whether the patient eventually showed up to the original appointment, canceled the appointment some time later, or missed the appointment. Canceled appointment slots are offered to new patients when possible.

\begin{figure}[h!]
   \centering
   \caption{Average number of arrivals in a week.}
   \label{fig:arrivals}
     \includegraphics[width=0.5\textwidth]{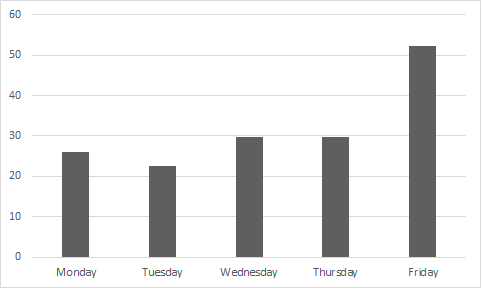}
\end{figure}

The average number of patients who arrive to make appointments on each day is shown in Figure \ref{fig:arrivals}. 
The actual arrival pattern is highly non-stationary, as the average number of arrivals on Friday is about twice that on other days. Our Marginal Allocation Algorithm gracefully handles this inherent non-stationarity.

We assume that there are two sessions on each day, a morning and an afternoon session. Each session on each day corresponds to a resource in our model. About $98\%$ appointments were scheduled into sessions on Monday through Thursday. We ignore the $2\%$ of appointments scheduled into other sessions because there are insufficient data to perform accurate analysis for these sessions. In other words, we set the capacity of sessions on Friday, Saturday and Sunday to be $0$. The capacity of sessions from Monday to Thursday are set based on the actual number of appointments made on these days, which is about 23 appointments per session. We will vary the capacity values in some of our experiments.

In this numerical experiment we do not model rescheduling, and treat each rescheduled appointment as an independent request. We also do not model the reuse of canceled appointment slots. Canceled slots are reused in practice, resulting in more efficient use of capacity.  In this way, our algorithms are at a disadvantage compared to actual practice because it has less capacity at its disposal.

We assume that the higher the probability that a patient will show up for a session, the more preferred the session is. Thus, we use {\it show probabilities} as a proxy for patient preferences for each session in a week. Specifically, we define the \emph{reward} of assigning a patient who arrives in period $i$ to a session $j$ as
\begin{equation} r_{ij} = \parbox{12cm}{Probability that the patient arriving in period $i$ will show up in session $j$ without canceling the appointment some time later or missing the appointment eventually.}\label{eq:CaseStudyBenefit} \end{equation}  
This definition of reward value does not capture all practical concerns, but it gives a good sense of scheduling effectiveness. The higher the measure is, the fewer no-shows and cancellations are likely to result, and the fewer appointments slots are potentially wasted.  In practice, operators try to subjectively assign appointments to accommodate patient preferences while maintaining a high level of utilization of capacity.  Because operators decisions are decentralized, they do not follow a precise and uniform procedure. However, our definition of reward is compatible with the goals of the actual system.

We estimate the show probabilities as a function of 3 factors: the day of the week, the time of day (morning/afternoon) and the number of days of wait starting from the patient's arrival to the actual appointment. In the first part of our experiment, we assume that patients have identical preferences in the sense that any two patients arriving on the same day will have the same reward values for each open session.  Thus, patients differ only in their time of arrival.

Both of the above assumptions regarding the homogeneity of preferences and the usefulness of show probabilities as indicators of preferences are strong assumptions.  We are aware that the show probabilities are imperfect substitute for actual preferences.  They also only express an average measure of preference.  A finer experiment would take into account actual preferences and variability of preferences among patients.  However, we believe that our experiment is still valuable in indicating the value of using online algorithms.  In a sense, our online algorithms are at a disadvantage compared to real practice because in practice, appointments were made taking into account actual preferences, whereas our online algorithms "know" only the show probabilities.

\begin{figure}[h!]
   \caption{Show probabilities of appointment slots assigned to patients who arrived on the previous Thursday.}
   \label{fig:showrate}
   \centering
     \includegraphics[width=0.5\textwidth]{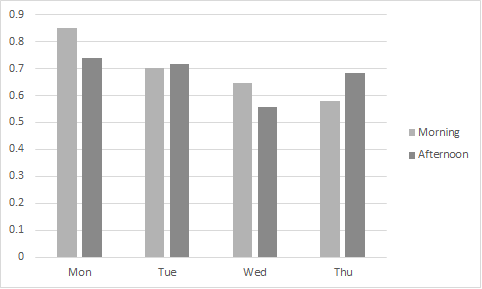}
\end{figure}

\begin{figure}[h!]
   \caption{Show probabilities as functions of number of days to wait before getting service.}
   \label{fig:showratelong}
   \centering
     \includegraphics[width=0.5\textwidth]{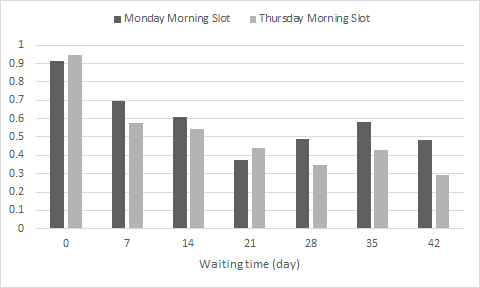}
\end{figure}

Figure \ref{fig:showrate} illustrates the show probabilities of patients who arrive on a Thursday to make appointments for the following week.
We can see that, in general, the shorter the wait is in days, the higher the show probability is.  Figure \ref{fig:showratelong} illustrates the show probability as a function of number of days to wait before getting service. The show probabilities range from as low as $27\%$, for appointments made more than two months into the future, to as high as $97\%$, for same-day visits.  
Table \ref{tab:showratetable} shows more show probabilities as a function of waiting time and day of week of the appointment.

\begin{table}[h!]
\caption{Show probabilities for morning sessions, as a function waiting time and day of week of the appointment. Some cells are NA because there is no patient arrival during weekends.}
\label{tab:showratetable}
\centering
\begin{tabular}{|c|c|c|c|c|c|c|c|c|c|}
\hline
 & \multicolumn{9}{|c|}{Number of days waiting}\\
\hline
Day of Week of Appointment&0 & 1 & 2 & 3 & 4 & 5 & 6 & 7 & 8\\
\hline
Mon	&$	91\%	$&	NA	&NA	&$	81\%	$&$	85\%	$&$	78\%	$&$	82\%	$&$	70\%	$&$	69\%	$\\
\hline
Tue	&$	78\%	$&$	83\%	$&NA	&NA	&$	62\%	$&$	70\%	$&$	73\%	$&$	58\%	$&$	53\%	$\\
\hline
Wed	&$	97\%	$&$	61\%	$&$	46\%	$&NA	&NA	&$	57\%	$&$	65\%	$&$	52\%	$&$	50\%	$\\
\hline
Thur	&$	95\%	$&$	67\%	$&$	41\%	$&$	50\%	$&NA	&NA	&$	60\%	$&$	58\%	$&$	57\%	$\\
\hline
\end{tabular}
\end{table}

      We used a 12-week period from March to May in 2013 as our time horizon.  An appointment reminder system was in use during this time. There are 2032 patients scheduled during this horizon according to our data. We use the sample consisting of these 2032 patient arrivals to simulate the performance of the following scheduling policies.

\begin{itemize}
 \item The Marginal Allocation Algorithm (MAA). The arrival rates, which are inputs of the algorithm, are estimated using our one-year data in 2013. The average number of arrivals in each day of week has been shown earlier in Figure \ref{fig:arrivals}. 
 \item The Marginal Allocation Algorithm with estimation error $\alpha\%$ (MAA-$\alpha\%$). This algorithm uses reward values (\ref{eq:CaseStudyBenefit}) that are each randomly and independently perturbed by relative errors drawn from a uniform distribution over $[-\alpha\%, \alpha\%]$. The total reward earned by this algorithm is computed using the unperturbed reward values.  We include these algorithms to test the impact of our parameter estimation errors on the performance comparison with actual practice.
\item The Separation Algorithm.
   \item  The outcome of actual practice used in hospitals. The total reward earned by the actual strategy is also calculated using the reward values defined in (\ref{eq:CaseStudyBenefit}).
   \item  The greedy policy, which always assigns a patient to the available session that is most preferred by the patient, as indicated by the show probability of the session. It captures a naive but easily implementable policy when a scheduler is aware of patient preferences. 
   \item The bid-price policy, which uses the optimal dual variables of LP (\ref{eq:LP2}) corresponding to the capacity constraints as the bid prices.  It assigns an arriving customer to the resource with the lowest price smaller than or equal to the revenue that the customer brings.  This heuristic is a widely used heuristic in resource-allocation problems.  
   \end{itemize}

In our first experiment, we do not consider overbooking and cancellations. The capacity of each session is set to be the number of appointments made in practice. In other words, we assume that the actual practice fully utilizes the capacity of all resources. Furthermore, we assume that patients arriving on the same day have homogeneous reward values.  

Since we use show probability as the reward of scheduling a patient, the total reward that a scheduling policy earns from the total $2032$ patients is equal to the  expected number of patients, among $2032$, who will show up to the original appointments. In particular, since the show probabilities are themselves estimated based on the scheduling of the actual practice, the total reward earned by the actual practice is just the actual number of patients, out of the total $2032$, that showed up during the horizon.

For each scheduling policy, we report as its performance the ratio of total reward to the total number $2032$ of arrivals. This ratio represents the overall percentage of patients who will show up. Table \ref{tab:numerical} summarizes the performance of all scheduling policies we consider.                                                                                      
We can see that our Marginal Allocation Algorithm performs the best, and in particular, gives more than $30\%$ improvement over the actual practice, according to our performance measure. It is noteworthy that the greedy and bid-price policies do not have performance guarantees and can perform arbitrarily badly. In contrast, our Marginal Allocation Algorithm has not only a provable performance guarantee, but also good empirical performance. 
      

      The strength of our Marginal Allocation Algorithm is more directly reflected in comparison with the greedy policy. The greedy policy can be carried out by anyone as long as the patient preferences are exploited. Our Marginal Allocation Algorithm, which does smart reservation, gives $12.9\%$ empirical improvement in scheduling efficiency over this heuristic.  Note that in this experiment, all patients have the same priority.  Our Marginal Allocation Algorithm is likely to exhibit much higher rewards when there are more patient types to consider because it can make more intelligent tradeoffs among the types than the greedy policy can. 	Remarkably, our Marginal Allocation Algorithm can be implemented as easily as the greedy policy. In the greedy approach, the scheduler has to be given a number representing estimated patient preference for each session. In our Marginal Allocation Algorithm, the scheduler also needs to be given only one number, namely the marginal value of reward function, for each session. 
      
\begin{table}
\begin{center}
\caption{The empirical performance of different scheduling policies.}
\label{tab:numerical}
\begin{tabular}{|c|c|}
\hline
Scheduling Policy & Performance of scheduling policies relative to LP upper bound\\
\hline
Actual Strategy & $67\%$\\
\hline
Greedy & $81\%$\\
\hline
Bid-Price Heuristic &  $89\%$\\
\hline
Separation Algorithm & $80\%$\\
\hline
MAA & $92\%$\\
\hline
MAA-$5\%$ & $91\%$\\
\hline
MAA-$10\%$ & $88\%$\\
\hline
MAA-$20\%$ & $83\%$\\
\hline
MAA-$40\%$ & $74\%$\\
\hline
\end{tabular}
\end{center}
\end{table}

\subsection{Consideration for Overbooking}
Starting from the numerical settings in the previous section, we study the practice of overbooking.  Let $A_j$ be the actual number of patients who are assigned to session $j$. We assume that the actual strategy overbooks each session by a constant ratio, and thereby treat $C_j = \alpha A_j$ as the actual capacity of session $j$, where $\alpha\in [0,1]$ is a scaling parameter that we vary in the numerical experiment.

We define the no-show probability as
\[\def\arraystretch{0.5}
 P_{NS} = \frac{\begin{array}{c}\text{Total number of no-shows} \\+ \\\text{ Total number of appointments that are canceled no more than}\\\text{ 2 days prior to the appointment time}\end{array}}{\text{Total number of appointments}}.\]

The number is $26.89\%$ as estimated from the data for Clinical Genetics.

A common practice is to take advantage of such high no-show probability by scheduling more patients to a session than its actual capacity can handle. Using terminology defined in Section \ref{sec:overbooking}, we use $P_{NS}$ as the no-show probability for every session. We also vary the no-show penalty $D$ in our experiments in the range $[2,10]$. In this way, the pair $(\alpha, D)$ tunes the cost (\ref{eq:overbookingcost}) of overbooking each session. The previous experiment corresponds to the case $\alpha = 1, D = \infty$.

Now the total reward of a scheduling policy is equal to the sum of all reward values  (\ref{eq:CaseStudyBenefit}), i.e., show probabilities, earned from patients less the overbooking costs (\ref{eq:overbookingcost}). In particular, we apply the function (\ref{eq:overbookingcost}) of overbooking cost to the actual practice as well. That is, in our experiment the total overbooking costs incurred under the actual practice does not depend on the actual overbooked number of patients, but rather on the expected costs (\ref{eq:overbookingcost}) estimated a priori. The performance of each scheduling policy is reported as its total reward relative to the total reward of the actual practice.

Table \ref{tab:numerical2} summarizes the performance of scheduling policies when $\alpha = 0.75$ and $D$ ranges from $2$ to $15$. Generally the performance of all policies decreases as the penalty $D$ increases because of the reduced reward of overbooking. 

Table \ref{tab:numerical3} reports the performance of scheduling policies when $D = 3$ and $\alpha$ increases from $70\%$ to $100\%$. The performance of all the scheduling policies reaches a limit for large values of $\alpha$. This is because when $\alpha$ is large, there is a large surplus of capacity associated with low overbooking costs. In such cases, scheduling policies virtually cannot see any capacity constraint, and thus have very good performance. Overall, for all values of $\alpha$, our Marginal Allocation Algorithm performs at least $30\%$ better than actual practice.
\begin{table}
\begin{center}
\caption{The total reward of scheduling policies relative to LP upper bound under different values of penalty $D$. $\alpha = 0.75$.}
\label{tab:numerical2}
\begin{tabular}{|c|c|c|c|c|c|}
\hline
$D$ & Actual Strategy & Greedy & Bid-Price Heuristic & Separation Alg. & MAA\\
\hline
2	&$	70.1\%	$&$	81.5\%	$&$	89.0\%	$&$	82.1\%	$&$	93.2\%	$\\
\hline
3	&$	68.7\%	$&$	80.6\%	$&$	86.7\%	$&$	82.0\%	$&$	92.4\%	$\\
\hline
4	&$	66.8\%	$&$	80.0\%	$&$	86.6\%	$&$	82.5\%	$&$	92.3\%	$\\
\hline
5	&$	64.5\%	$&$	79.5\%	$&$	87.2\%	$&$	82.8\%	$&$	92.3\%	$\\
\hline
6	&$	62.2\%	$&$	79.2\%	$&$	88.7\%	$&$	82.5\%	$&$	92.0\%	$\\
\hline
7	&$	59.7\%	$&$	78.9\%	$&$	88.0\%	$&$	82.6\%	$&$	92.0\%	$\\
\hline
8	&$	57.1\%	$&$	78.7\%	$&$	88.6\%	$&$	82.5\%	$&$	92.0\%	$\\
\hline
9	&$	54.5\%	$&$	78.2\%	$&$	88.4\%	$&$	82.4\%	$&$	92.0\%	$\\
\hline
10	&$	51.8\%	$&$	77.8\%	$&$	88.4\%	$&$	82.1\%	$&$	91.5\%	$\\
\hline
\end{tabular}
\end{center}
\end{table}

\begin{table}
\begin{center}
\caption{The total reward of scheduling policies relative to LP upper bound under different values of $\alpha$. $D = 3$.}
\label{tab:numerical3}
\begin{tabular}{|c|c|c|c|c|c|}
\hline
$\alpha$  & Actual Strategy & Greedy & Bid-Price Heuristic & Separation Alg. & MAA\\
\hline
$	70\%	$&$	62.7\%	$&$	77.2\%	$&$	88.3\%	$&$	82.9\%	$&$	92.2\%	$\\
\hline
$	75\%	$&$	68.7\%	$&$	80.6\%	$&$	86.7\%	$&$	82.0\%	$&$	92.4\%	$\\
\hline
$	80\%	$&$	70.8\%	$&$	83.9\%	$&$	88.9\%	$&$	81.8\%	$&$	93.4\%	$\\
\hline
$	85\%	$&$	71.4\%	$&$	88.9\%	$&$	92.5\%	$&$	82.6\%	$&$	94.5\%	$\\
\hline
$	90\%	$&$	71.1\%	$&$	91.4\%	$&$	94.4\%	$&$	83.3\%	$&$	94.8\%	$\\
\hline
$	95\%	$&$	70.7\%	$&$	92.5\%	$&$	95.9\%	$&$	84.6\%	$&$	95.8\%	$\\
\hline
$	100\%	$&$	70.5\%	$&$	93.2\%	$&$	95.7\%	$&$	85.7\%	$&$	96.2\%	$\\
\hline
\end{tabular}
\end{center}
\end{table}

\subsection{Consideration for Patient Availability}

In the previous numerical experiments, patients who arrive in the same periods are treated as identical. However, in reality there is variability among patients' availability. In this section, we capture this variability by simulating a particular chosen patient's availability for a particular session of the week as being drawn from a given distribution. This experiment tests whether more complex heterogeneous patient types affect the comparative performance of our algorithm.

We model the heterogeneity of patient availability as follows. A patient cannot be assigned to a session if he is unavailable for it. Otherwise, the reward for the session is still the show probability as modeled in the previous sections. We assume that each patient has the same availability pattern for every week. A patient is available for any session with probability $P_A$, and this event is independent of the availability for other sessions in the same week. We vary $P_A$ from $15\%$ to $100\%$ to test the performance of all the scheduling policies we consider. When $P_A= 100\%$, the problem is reduced to the one in the last section, in which a patient can be assigned to any session.

Since we model 8 sessions in a week, one in the morning and one in the afternoon from Monday to Thursday (recall that there were very few appointments scheduled for Friday), each patient's availability can be represented by an 8-dimension binary vector. Then, patients arriving in each period are further divided into $2^8$ patient types, with $r_{i,k,j}=0$ if a patient of type $k \in \{1,2,...,2^8\}$ arriving in period $i$ is not available for session $j$.

We assume that the sessions offered by actual practice to patients were all available, so that the total reward of actual practice is not affected by this newly modeled feature. The performance of each of the remaining scheduling policies is the averaged total reward over $10,000$ runs of simulation. In each simulation we draw the same $2032$ number of arrivals from data, but we randomly generate patient availability.
For $P_A$ ranging from $15\%$ to $100\%$, Table \ref{tab:numerical4} shows the performance of scheduling policies relative to the performance of actual practice.  The relative performance is better for higher values of $P_A$, as there is more flexibility in scheduling when patients are available to more sessions. Even when $P_A$ is as small as $15\%$, our Marginal Allocation Algorithm still performs $8\%$ better than actual practice. The gap gradually increases to more than $40\%$ as $P_A$ increases.

\begin{table}
\begin{center}
\caption{The total reward of scheduling policies relative to LP upper bound under different values of $P_A$. $D = 3$, $\alpha = 0.7$.}
\label{tab:numerical4}
\begin{tabular}{|c|c|c|c|c|c|}
\hline
$P_A$  & Actual Strategy & Greedy & Bid-Price Heuristic & Separation Alg. & MAA\\
\hline
$	15.00\%	$&$	88.8\%	$&$	89.6\%	$&$	96.0\%	$&$	93.1\%	$&$	96.4\%	$\\
\hline												
$	20.00\%	$&$	77.6\%	$&$	86.0\%	$&$	93.5\%	$&$	90.9\%	$&$	95.1\%	$\\
\hline												
$	25.00\%	$&$	72.5\%	$&$	83.4\%	$&$	92.1\%	$&$	89.5\%	$&$	94.2\%	$\\
\hline												
$	30.00\%	$&$	70.0\%	$&$	81.5\%	$&$	91.4\%	$&$	88.9\%	$&$	93.7\%	$\\
\hline												
$	35.00\%	$&$	68.4\%	$&$	80.2\%	$&$	91.4\%	$&$	88.5\%	$&$	93.5\%	$\\
\hline												
$	40.00\%	$&$	67.2\%	$&$	79.3\%	$&$	90.8\%	$&$	87.9\%	$&$	93.2\%	$\\
\hline												
$	45.00\%	$&$	66.3\%	$&$	78.7\%	$&$	91.5\%	$&$	87.6\%	$&$	93.0\%	$\\
\hline												
$	50.00\%	$&$	65.7\%	$&$	78.2\%	$&$	90.7\%	$&$	86.8\%	$&$	92.7\%	$\\
\hline												
$	55.00\%	$&$	65.1\%	$&$	77.8\%	$&$	90.5\%	$&$	86.0\%	$&$	92.7\%	$\\
\hline												
$	60.00\%	$&$	64.7\%	$&$	77.6\%	$&$	89.3\%	$&$	85.5\%	$&$	92.7\%	$\\
\hline												
$	65.00\%	$&$	64.3\%	$&$	77.5\%	$&$	88.7\%	$&$	85.0\%	$&$	92.8\%	$\\
\hline												
$	70.00\%	$&$	64.0\%	$&$	77.5\%	$&$	88.9\%	$&$	84.5\%	$&$	92.5\%	$\\
\hline												
$	75.00\%	$&$	63.7\%	$&$	77.6\%	$&$	88.6\%	$&$	84.0\%	$&$	92.4\%	$\\
\hline												
$	80.00\%	$&$	63.4\%	$&$	77.6\%	$&$	87.7\%	$&$	83.7\%	$&$	92.3\%	$\\
\hline												
$	85.00\%	$&$	63.2\%	$&$	77.6\%	$&$	88.3\%	$&$	83.5\%	$&$	92.4\%	$\\
\hline												
$	90.00\%	$&$	63.0\%	$&$	77.5\%	$&$	88.4\%	$&$	83.3\%	$&$	92.5\%	$\\
\hline												
$	95.00\%	$&$	62.9\%	$&$	77.5\%	$&$	88.0\%	$&$	83.1\%	$&$	92.4\%	$\\
\hline												
$	100.00\%	$&$	62.7\%	$&$	77.2\%	$&$	88.3\%	$&$	82.9\%	$&$	92.2\%	$\\
\hline
\end{tabular}
\end{center}
\end{table}

\section{Appendix: Omitted Proofs}

\noindent {\bf Proof of Theorem \ref{thm:nboundthm1}.}
\proof{Proof.}
Note that the distribution of $\{R(t)\}_{t\geq 0}$ is determined by the time points $t_1,t_2,...,t_{k-1}$. In particular, for $t\in(t_{i},t_{i+1})$, the value of $\bP(R(t) = i)$ is only determined by $t_1,t_2,...,t_i$.

Given any value $\beta \in (0,1)$, we can construct a sequence of those time points $t_1, t_2,...,t_{k-1}$ recursively based on the following condition
\begin{equation}\label{eq:nboundpft1}  
\int_{t_{i}}^{t_{i+1}} \bP(R(t) = i) dt =  \frac{1}{\beta} - 1, \,\,\forall i = 0,1,...,k-2.
\end{equation}
Here $t_{i}$ is when the barrier is increased to position $i$, and is thus the first time that $\bP(R(t) = i)$ becomes positive. Given $t_1,t_2,...,t_i$, this condition sets the value for $t_{i+1} = t_{i+1}(\beta)$ by requiring that the area under the function $\bP(R(t) = i)$ for $t \in [t_i, t_{i+1}]$ is exactly $1/\beta - 1$.

According to the above construction, since $\bP(R(t) = i)$ is a continuous function of $t$, the time points $t_1,t_2,..., t_{k-1}$ must change continuously in $\beta$.

 Furthermore, when $\beta \to 1$, i.e., the area under the function $\bP(R(t) = i)$ for $t \in [t_i, t_{i+1}]$ tends to $0$ for each $i = 0,1,...,k-2$, we must have $t_{i+1} - t_i \to 0$ for each $i = 0,1,...,i-2$. This implies that $t_{k-1} \to t_0 = 0$. On the other hand, when $\beta \to 0$, we have $1/{\beta} -1 \to \infty$, so the area under $\bP(R(t)=i)$ for $t\in[t_i,t_{i+1}]$ can be arbitrarily large. In other words, by tuning the value of $\beta$, we can set $t_{k-1}$ to be any value within $(0,k)$.

Therefore, there must exist some $\beta \in (0,1)$ such that $t_{k-1}$ satisfies
\[ \int_{t_{k-1}}^{t_{k}} \bP(R(t) = k-1) dt =  \frac{1}{\beta} - 1.\]
Let $\beta^*$ be such a value that satisfies this condition. Set $\alpha^*_i(t) =  \bP(R(t) = i) \beta^*$. We next prove that this construction of $\beta^*$ and $\alpha^*_i(t)$, for $i=0,1,...,k-1$ and $t\in[0,k]$, satisfies the constraints of (\ref{eq:dual}).

First of all, for $t \leq t_1$, we have
\begin{align*}
 \alpha_0^*(t) + \int_0^t \alpha_0^*(s) ds & = \beta^* \bP(R(t) = 0) + \int_0^t \beta^* \bP(R(s) = 0) ds\\
 & = \beta^* \cdot 1 + \beta^* \int_0^t \bP(R(s) = 0) ds\\
 & \leq \beta^* + \beta^* \int_0^{t_1} \bP(R(s) = 0) ds\\
 & = \beta^* + \beta^* ( 1 / \beta^* -1)\\
 & = 1.
  \end{align*}
Note that the inequality is tight when $t = t_1$. For $t > t_1$, since the barrier is above position $0$, the random process $R(t)$ is changing from state $R(t)=0$ to state $R(t)=1$ at rate $1$ (the transition happens when $N(t)$ increases by $1$). Thus, we must have, for $t > t_1$,
\begin{align*}
 &\frac{\partial \bP(R(t) = 0)}{\partial t} = - \bP(R(t) = 0)  \\
 \Longrightarrow & \bP(R(t) =0) - \bP(R(t_1) = 0) = -\int_{t_1}^t \bP(R(s) = 0) ds\\
 \Longrightarrow & \bP(R(t) = 0) + \int_0^t \bP(R(s) = 0) ds = \bP(R(t_1) = 0) + \int_0^{t_1} \bP(R(s) = 0) ds \\
 \Longrightarrow & \alpha^*_0(t) + \int_0^t \alpha_0^*(s) ds = \alpha^*_0(t_1) + \int_0^{t_1} \alpha^*_0(s)ds = 1.
 \end{align*}
Therefore, the first constraint of (\ref{eq:dual}) holds and is tight for $t\geq t_1$.

To prove that the second constraint also holds, we recursively look at $i=1,2,...,k-1$. Recall that $t_i$ is the first time that $\bP(R(t) = i)$ becomes positive. Thus for $t \leq t_i$ we have $\bP(R(t) = i) = 0$ and
\[ \alpha^*_i(t) + \int_0^t \alpha_i^*(s) ds = \beta^* \bP(R(t) = i) + \int_0^t \beta^* \bP(R(s) = i) ds = 0.\]

For $t \in [t_{i}, t_{i+1}]$, we have
\begin{align*}
\alpha_i^*(t) + \int_0^t \alpha_i^*(s) ds & = \beta^* \bP(R(t) = i) + \int_0^t \beta^* \bP(R(s) = i) ds \\
& = \beta^* \bP(R(t) = i) + \beta^* \int_{t_i}^{t} \bP(R(s) = i) ds\\
& \leq \beta^* \bP(R(t) = i) + \beta^*\int_{t_i}^{t_{i+1}} \bP(R(s) = i) ds\\
& = \beta^* \bP(R(t) = i) + \beta^* (1 / \beta^* - 1)\\
& = \beta^* \bP(R(t) = i) + \beta^* \int_{t_{i-1}}^{t_{i}} \bP(R(s) = i-1) ds.
\end{align*}

When $t \in [t_i, t_{i+1}]$ and $R(t) = i$, the random process is actively bounded above by the barrier, so the probability $\bP(R(t) =i)$, as a function of $t$, can only increase due to the transition from state $R(t) = i-1$ to $R(t) = i$. The rate of this transition is $1$. Thus, we have $\bP(R(t) = i) = \int_{t_i}^t \bP(R(s) = i-1) ds$, which leads to
\begin{align}
\begin{split}
 \alpha_i^*(t) + \int_0^t \alpha_i^*(s) ds \leq & \beta^* \bP(R(t) = i) + \beta^* \int_{t_{i-1}}^{t_{i}} \bP(R(s) = i-1) ds\\
 = &\beta^* \int_{t_i}^t \bP(R(s) = i-1) ds + \beta^* \int_{t_{i-1}}^{t_{i}} \bP(R(s) = i-1) ds\\
= & \beta^* \int_0^t \bP(R(s) = i-1) ds\\
= & \int_0^t \alpha_{i-1}^*(s) ds.
\end{split}
 \end{align}
Note that the inequality is tight when $t = t_{i+1}$.

For $t > t_{i+1}$, the barrier is above $i$, so the random process $R(t)$, if still at state $R(t) = i$, is not actively bounded by the barrier. Thus  the state $R(t) = i$ is involved in two transitions: from state $i$ to $i+1$, and from $i-1$ to $i$. More precisely, we have for $t > t_{i+1}$,
\[ \frac{\partial \bP(R(t) = i)}{\partial t} = \bP(R(t) = i-1) - \bP(R(t) = i)\]
\begin{align*}
 \Longrightarrow & \bP(R(t) = i) + \int_0^t \bP(R(s) = i) ds - \int_0^t \bP(R(s) = i-1) ds\\
& = 0.
\end{align*}
\begin{equation}\label{eq:nboundpfbbb}
 \Longrightarrow \alpha_i^*(t) + \int_0^t \alpha_i^*(s) ds = \int_0^t \alpha_{i-1}^*(s) ds.
 \end{equation}
This proves that the second constraint of (\ref{eq:dual}) holds (and is tight for $t\geq t_{i+1}$, for each $i=1,2,...,k-1$, respectively). Finally, the last constraint of (\ref{eq:dual}) trivially holds because $\sum_{i=0}^{k-1} \bP(R(t) = i) = 1 \Longrightarrow \beta^* = \sum_{i=0}^{k-1} \alpha^*_i(t)$.

To prove (\ref{eq:nboundthm1eq1}), we can deduce that
\begin{align*}
& \beta^* \sum_{i=0}^{k-1} i \bP(R(k) = i)\\
= & \sum_{i=0}^{k-1} i \alpha_i^*(k)\\
= & \sum_{i=0}^{k-1} i \left[ \int_0^k \alpha^*_{i-1}(s) ds - \int_0^k \alpha^*_i(s)ds\right]\,\,\,\text{(by (\ref{eq:nboundpfbbb}))}\\
= &  \sum_{i=0}^{k-1} \int_0^k \alpha^*_i(s) ds - k \int_0^k \alpha^*_{k-1}(s) ds\,\,\,\text{(by canceling identical terms)} \\
= & \int_0^k \left( \sum_{i=0}^{k-1} \alpha^*_i(s) \right) ds - k \beta^* \int_0^k \bP(R(s) = k-1) ds\\
= &  \int_0^k \beta^* ds - k\beta^* (1 / \beta^* - 1)\\
= & 2 k\beta^* - k.
\end{align*}
We can then easily obtain (\ref{eq:nboundthm1eq1}) by rearranging terms.
\halmos\endproof

\noindent {\bf Proof of Lemma \ref{lm:nboundlm4}.}
\proof{Proof.}

\begin{align*}
& \mathbf{E}[I_l| R(l) = l-j]\\
= & \int_l^{l+1} \bP(R(s) = l | R(l) = l-j) ds\\
= & \int_0^{1} P_{\geq j}(s) ds.
\end{align*}

Similarly, $\mathbf{E}[I_l | R(l) = l-j-1] = \int_0^1 P_{\geq j+1}(s) ds$. Thus,
\begin{align*}
& \mathbf{E}[I_l | R(l) = l-j] - \mathbf{E}[I_l | R(l) = l-j-1] \\
= & \int_0^{1} P_{\geq j}(s) ds - \int_0^1 P_{\geq j+1}(s) ds\\
= & \int_0^1 P_{j}(s) ds\\
= & \int_0^1 e^{-s} \frac{s^j}{j!} ds\\
= & \sum_{\nu = j+1}^\infty e^{-1} \frac{1}{\nu!}\\
= & P_{\geq j+1}(1).
\end{align*}
\halmos\endproof

\noindent {\bf Proof of Lemma \ref{lm:nboundlmkeybound}.}
\proof{Proof.}
Fix any $i \in \{l,l+1,...,k-1\}$. For ease of notation, define
\[ \Delta_{d,j} \equiv \mathbf{E}[I_i | R(d) = d-j] - \mathbf{E}[I_{i}|R(d) = d-j-1] \]
to be the increment in the expected time that $R(t)$ stays at the barrier during $[t_i, t_{i+1}]$, when the state at time $t = d$ changes from $R(d) = d-j-1$ to $R(d) = d-j$, for all $ d= l,l+1,...,i$ and $j = 0,1,...,d-1$. 

From Lemma \ref{lm:nboundlm4} we know that $\Delta_{i,j} = P_{\geq j+1}(1)$. Furthermore, for $d<i$ and $d \geq l$, the value of $\mathbf{E}[I_i | R(d) = d-j] $ can be recursively computed by conditioning on $R(d+1)$, i.e., on the movement of the random process during time $[d, d+1]$. Precisely,
\[ \mathbf{E}[I_i | R(d) = d-j] = \sum_{\nu=0}^j P_\nu(1) \mathbf{E}[I_i |  R(d+1) = d-j+\nu] + \sum_{\nu = j+1}^\infty P_\nu(1) \mathbf{E}[I_i | R(d+1) = d].  \]
Here, for example, $R(d+1) = d-j+\nu$ represents the condition where the random process $R(t)$ moves $\nu$ steps upwards during time $[d,d+1]$; $R(d+1) = d$ represents the condition where the barrier is active at time $t = d+1$.

Similarly, 
\begin{align*}
 \mathbf{E}[I_i | R(d) = d-j-1] &= \sum_{\nu=0}^{j+1} P_\nu(1) \mathbf{E}[I_i |  R(d+1) = d-j-1+\nu] + \sum_{\nu = j+2}^\infty P_\nu(1) \mathbf{E}[I_i | R(d+1) = d]  \\
 & = \sum_{\nu=0}^{j} P_\nu(1) \mathbf{E}[I_i |  R(d+1) = d-j-1+\nu] + \sum_{\nu = j+1}^\infty P_\nu(1) \mathbf{E}[I_i | R(d+1) = d] .
 \end{align*}
The above two equations lead to the following recursion for $\Delta_{d,j}$
\begin{align}
\begin{split}\label{eq:nboundrecursiondelta}
 \Delta_{d,j}& = \mathbf{E}[I_i | R(d) = d-j] - \mathbf{E}[I_i | R(d) = d-j-1] \\
 & =  \sum_{\nu = 0}^j P_\nu(1) [  \mathbf{E}[I_i |  R(d+1) = d-j+\nu] -\mathbf{E}[I_i |  R(d+1) = d-j-1+\nu]]\\
 & =  \sum_{\nu = 0}^j P_\nu(1) \Delta_{d+1,j-\nu+1}.
 \end{split}
 \end{align}

Note that in order to prove the lemma, we need to show $\Delta_{l,0} / \Delta_{l,1} \geq 1/e$. To this end, we prove a stronger result by constructing a bound on $\Delta_{d,j} / \Delta_{d,j+1}$ for all $d = l,l+1,...,i$, and $j=0,1,...,d$. 
We construct the bounds using a sequence of `stationary' values $\Delta_{*,0}, \Delta_{*,1},...$, which are defined based on the recursion (\ref{eq:nboundrecursiondelta}) and are independent of $d$:
\[ \Delta_{*,0} = 1; \]
\begin{align}
\Delta_{*,j} = &\sum_{\nu = 0}^j P_\nu(1) \Delta_{*,j-\nu+1},\,\,\,\forall j = 0,1,2,... \label{eq:nboundpfdeltastarleft}\\
& \Longrightarrow  \left\{ \begin{array}{ll}\Delta_{*,1} = e \Delta_{*,0},\\ \Delta_{*,j+1} =(e-1) \Delta_{*,j} - \sum_{\nu = 2}^j \frac{1}{\nu!}\Delta_{*,j-\nu+1}, \,\,\forall j \geq 1.\end{array} \right. \label{eq:nboundpfdeltastar}
\end{align}

We next prove that
\begin{equation}\label{eq:nboundpfdelta}
 \frac{\Delta_{d,j}}{\Delta_{d,j+1}} \geq \frac{\Delta_{*,j}}{\Delta_{*,j+1}}
 \end{equation}
by induction on $d$.
\begin{itemize}
\item
First, we prove that (\ref{eq:nboundpfdelta}) holds for $d=i$, by showing that $\Delta_{i,j}$ is decreasing in $j$ but $\Delta_{*,j}$ is increasing in $j$. 

By Lemma \ref{lm:nboundlm4}, $\Delta_{i,j} = P_{\geq j+1}(1) > P_{\geq j+2}(1) = \Delta_{i,j+1}$, which means $\Delta_{i,j}$ is decreasing in $j$.  

From (\ref{eq:nboundpfdeltastar}) we know that $\Delta_{*,0} = e^{-1} \Delta_{*,1} < \Delta_{*,1}$. Provided that $\Delta_{*,\nu} \leq \Delta_{*,j}$ for all $\nu \leq j$ and some $j\geq 1$, we can deduce from (\ref{eq:nboundpfdeltastar}),
\[\frac{\Delta_{*,j+1}}{\Delta_{*,j}} \geq e-1  - \sum_{\nu = 2}^j \frac{1}{\nu!} \geq e-1 - \sum_{\nu=2}^\infty \frac{1}{\nu!} = e-1-(e-2)= 1.\]
Thus, $\Delta_{*,j}$ is increasing in $j$, which finishes the proof that (\ref{eq:nboundpfdelta}) holds when $d=i$.
\item When $d < i$,
\begin{align}
\nonumber& \Delta_{d,j} \Delta_{*,j+1} - \Delta_{d,j+1} \Delta_{*,j}\\
\nonumber= &  \left(\sum_{\nu_1 = 0}^j P_{\nu_1}(1) \Delta_{d+1,j-\nu_1+1}\right)\left( \sum_{\nu_2 = 0}^{j+1} P_{\nu_2}(1) \Delta_{*,j-\nu_2+2}\right)\\
\nonumber& \,\,\,\,\,\,\,\,\,\, - \left(\sum_{\nu_1 = 0}^j P_{\nu_1}(1) \Delta_{*,j-\nu_1+1}\right)\left( \sum_{\nu_2 = 0}^{j+1} P_{\nu_2}(1) \Delta_{d+1,j-\nu_2+2}\right) \,\,\text{(by (\ref{eq:nboundrecursiondelta}) and (\ref{eq:nboundpfdeltastarleft}))} \\
\nonumber= &  \left(\sum_{\nu_1 = 0}^j P_{\nu_1}(1) \Delta_{d+1,j-\nu_1+1}\right) P_0(1) \Delta_{*,j+2} - \left(\sum_{\nu_1 = 0}^j P_{\nu_1}(1) \Delta_{*,j-\nu_1+1}\right) P_0(1) \Delta_{d+1,j+2}\\
\nonumber& \,\,\,\,\, + \left(\sum_{\nu_1 = 0}^j P_{\nu_1}(1) \Delta_{d+1,j-\nu_1+1}\right)\left( \sum_{\nu_2 = 0}^{j} P_{\nu_2+1}(1) \Delta_{*,j-\nu_2+1}\right)\\
\nonumber& \,\,\,\,\, - \left(\sum_{\nu_1 = 0}^j P_{\nu_1}(1) \Delta_{*,j-\nu_1+1}\right)\left( \sum_{\nu_2 = 0}^{j} P_{\nu_2+1}(1) \Delta_{d+1,j-\nu_2+1}\right)\\
& = \sum_{\nu_1 = 0}^j P_{\nu_1}(1) P_0(1) \left(   \Delta_{d+1,j-\nu_1+1} \Delta_{*,j+2} - \Delta_{*,j-\nu_1+1} \Delta_{d+1,j+2} \right) \nonumber\\ 
& \,\,\,\, + \sum_{\nu_1 = 0}^j \sum_{\nu_2 = 0} ^ {\nu_1-1} P_{\nu_1}(1) P_{\nu_2+1}(1)  \left(\Delta_{d+1,j-\nu_1 + 1} \Delta_{*,j-\nu_2+1} - \Delta_{*,j-\nu_1 + 1} \Delta_{d+1,j-\nu_2+1}\right)\nonumber\\
& \,\,\,\, +  \sum_{\nu_2 = 0}^j \sum_{\nu_1 = 0}^{\nu_2-1}  P_{\nu_1}(1) P_{\nu_2+1}(1)  \left(\Delta_{d+1,j-\nu_1 + 1} \Delta_{*,j-\nu_2+1} - \Delta_{*,j-\nu_1 + 1} \Delta_{d+1,j-\nu_2+1}\right) \nonumber \\
& \,\,\,\, + \sum_{\nu_1 = 0}^ j P_{\nu_1}(1) P_{\nu_1+1}(1)  \left(\Delta_{d+1,j-\nu_1 + 1} \Delta_{*,j-\nu_1+1} - \Delta_{*,j-\nu_1 + 1} \Delta_{d+1,j-\nu_1+1}\right) \nonumber \\
& = \sum_{\nu_1 = 0}^j P_{\nu_1}(1) P_0(1) \left(   \Delta_{d+1,j-\nu_1+1} \Delta_{*,j+2} - \Delta_{*,j-\nu_1+1} \Delta_{d+1,j+2} \right) \nonumber\\ 
& \,\,\,\, + \sum_{\nu_1 = 0}^j \sum_{\nu_2 = 0} ^ {\nu_1-1} P_{\nu_1}(1) P_{\nu_2+1}(1)  \left(\Delta_{d+1,j-\nu_1 + 1} \Delta_{*,j-\nu_2+1} - \Delta_{*,j-\nu_1 + 1} \Delta_{d+1,j-\nu_2+1}\right)\nonumber\\
& \,\,\,\, +  \sum_{\nu_2 = 0}^j \sum_{\nu_1 = 0}^{\nu_2-1}  P_{\nu_1}(1) P_{\nu_2+1}(1)  \left(\Delta_{d+1,j-\nu_1 + 1} \Delta_{*,j-\nu_2+1} - \Delta_{*,j-\nu_1 + 1} \Delta_{d+1,j-\nu_2+1}\right) \nonumber \\
& = \sum_{\nu_1 = 0}^j P_{\nu_1}(1) P_0(1) \left( \Delta_{d+1,j-\nu_1+1} \Delta_{*,j+2} - \Delta_{*,j-\nu_1+1} \Delta_{d+1,j+2} \right) \label{eq:nboundpfinduction1}\\
& \,\,\,\, + \sum_{\nu_1 = 0}^j \sum_{\nu_2 = 0} ^ {\nu_1-1}  \left(P_{\nu_1}(1) P_{\nu_2+1}(1) -P_{\nu_2}(1) P_{\nu_1+1}(1)\right)   \left(\Delta_{d+1,j-\nu_1 + 1} \Delta_{*,j-\nu_2+1} - \Delta_{*,j-\nu_1 + 1} \Delta_{d+1,j-\nu_2+1}\right) \label{eq:nboundpfinduction2}
\end{align}
Now using induction on $d+1$, we know that for $\nu_1 \geq 0$, 
\[\frac{\Delta_{d+1,j-\nu_1+1}}{\Delta_{d+1,j+2}} \geq \frac{\Delta_{*,j-\nu_1+1}}{\Delta_{*,j+2}},\]
and thus (\ref{eq:nboundpfinduction1})$\geq 0$. In  (\ref{eq:nboundpfinduction2}), since $\nu_1 > \nu_2$, we can use induction on $d+1$ again to obtain
\[ \frac{\Delta_{d+1,j-\nu_1+1}}{\Delta_{d+1,j-\nu_2+1}} \geq \frac{\Delta_{*,j-\nu_1+1}}{\Delta_{*,j-\nu_2+1}}\]
\[ \Longrightarrow \Delta_{d+1,j-\nu_1 + 1} \Delta_{*,j-\nu_2+1} - \Delta_{*,j-\nu_1 + 1} \Delta_{d+1,j-\nu_2+1} \geq 0.\]
Furthermore, since $\nu_1 > \nu_2$,
\[ P_{\nu_1}(1) P_{\nu_2+1}(1) -P_{\nu_2}(1) P_{\nu_1+1}(1) = P_{\nu_1}(1) P_{\nu_2}(1)\left(\frac{1}{\nu_2+1} - \frac{1}{\nu_1+1}\right)\geq 0.\]
Thus, (\ref{eq:nboundpfinduction2})$\geq 0$. In sum, we have shown $\Delta_{d,j} \Delta_{*,j+1} - \Delta_{d,j+1} \Delta_{*,j} \geq 0$, which finishes the induction proof for condition (\ref{eq:nboundpfdelta}).
\end{itemize}

This result (\ref{eq:nboundpfdelta}) directly leads to the statement of the Lemma
\[ \Delta_{l,0} \geq \Delta_{l,1} \frac{\Delta_{*,0}}{\Delta_{*,1}} = \Delta_{l,1} / e.\]

\halmos\endproof

\noindent {\bf Proof of Lemma \ref{lm:nboundlm5}.}

\proof{Proof.}
Fix any $i \in \{l,l+1,...,k-2\}$. By symmetry, we have $\mathbf{E}[I_i | R(l) = l] = \mathbf{E}[I_{i+1} | R(l+1) = l+1]$ and $\mathbf{E}[I_i |  R(l) = l-1] = \mathbf{E}[I_{i+1} | R(l+1) = l]$.

Thus,
\begin{align*}
& \mathbf{E}[I_i | R(l) = l] - \mathbf{E}[I_{i+1}|R(l) = l] - (\mathbf{E}[I_i | R(l) = l-1] - \mathbf{E}[I_{i+1} | R(l) = l-1])\\
= & \mathbf{E}[I_{i+1} | R(l+1) = l+1] -  \mathbf{E}[I_{i+1} | R(l+1) = l] - ( \mathbf{E}[I_{i+1}|R(l) = l] - \mathbf{E}[I_{i+1} | R(l) = l-1]))\\
= & \mathbf{E}[I_{i+1} | R(l+1) = l+1] -  \mathbf{E}[I_{i+1} | R(l+1) = l] - e^{-1}( \mathbf{E}[I_{i+1}|R(l+1) = l] - \mathbf{E}[I_{i+1} | R(l+1) = l-1]))\\
\geq & 0,
\end{align*}
where the last inequality follows from Lemma \ref{lm:nboundlmkeybound}.
\halmos\endproof

\noindent {\bf Proof of Lemma \ref{lm:nboundlmkey}.}

\proof{Proof.}
Given any sequence of times points $\bar t_1 \leq \bar t_2\leq \cdots \leq\bar t_l$ such that $\bar t_j \leq j$, $\forall j = 1,2,...,l$, we want to prove the lemma when $t_j = \bar t _j$ for $j \leq l$ and $t_j = j$ for $j > l$.

Fix any $i \in \{l,l+1,...,k-2\}$. Initially, set $t_j = j$ for all $j=1,2,...,k-1$, and we know that $\mathbf{E}[I_i] \geq \mathbf{E}[I_{i+1}]$ by Lemma \ref{lm:nboundlmmonotone}. We next prove that $\mathbf{E}[I_i] \geq \mathbf{E}[I_{i+1}]$ always holds when we reduce the value of $t_j$ from $j$ to $\bar t_j$ \emph{sequentially} for $j=1,2,...,l$.

Define  
 \[ I_i' = \left\{ \begin{array}{ll} I_i, & \text{ if } i > l \\ \int_l^{l+1} \mathbf{1}(R(s) = l) ds, & \text{ if } i = l. \end{array} \right.\]
By this definition, $I_i$ is different from $I_i'$ only if $i = l$ and $t_l < l$ (we create the definition of $I_i'$ so as to facilitate the proof for the case of $i=l$). Note that we always have $I_i \geq I_i'$.

Consider the result of reducing $t_j$ from $j$ to $\bar t_j$, when $t_d = d$ for all $d = j+1,j+2,...,k-1$. Let $\hat R(t), \hat I_i$ and $\hat I_i'$ be the values of $R(t), I_i$ and $I_i'$, respectively, before reducing $t_j$, i.e., when $t_j = j$. Let $\bar R(t), \bar I_i$ and $\bar I_i'$ be the new values of $R(t), I_i$ and $I_i'$, respectively, after reducing $t_j$ to $\bar t_j$.

Suppose $\mathbf{E}[\hat I_i] \geq \mathbf{E}[\hat I_{i+1}]$ holds. We next prove that $\mathbf{E}[\bar I_i] \geq \mathbf{E}[\bar I_{i+1}]$ also holds.

\begin{itemize}
\item  When $t_j = j$, since $j \leq l$, we must have $t_l = l$ and thus $\mathbf{E}[\hat I_i] = \mathbf{E}[\hat I_i'] $. 
\item We must have $\bP(\bar R(j) = \nu) = \bP(\hat R(j) = \nu)$ for all $\nu\leq j-2$, because if $\hat R(j) \leq j-2$, the random process does not touch the barrier during $[t_{j-1}, j]$, and is not affected by the change in $t_j$.
\item We must have $\bP(\bar R(j) = j-1) \leq \bP(\hat R(j) = j-1)$ and $\bP(\bar R(j) = j) \geq \bP(\hat R(j) = j)$ because when $t_j$ becomes smaller, there is more time for the random process $R(t)$ to jump from state $j-1$ up to $j$. Moreover,
\begin{equation}\label{eq:nboundlmpfchangeprob}
 \bP(\bar R(j) = j) - \bP(\hat R(j) = j) = \bP(\hat R(j) = j-1) - \bP(\bar R(j) = j-1) \geq 0.
 \end{equation}
\end{itemize}

Based on the above results, we can then deduce that ($\mathbf{E}[\bar I_i]$ is defined as $\mathbf{E}[I_i]$ given $t_j = \bar t_j$. Similarly, $\mathbf{E}[\hat I_i]$ is defined as $\mathbf{E}[I_i]$ given $t_j = j$)
\begin{align*}
& \mathbf{E}[\bar I_i] - \mathbf{E}[\bar I_{i+1}] \\
\geq &\mathbf{E}[\bar I_i'] - \mathbf{E}[\bar I_{i+1}] \\
= & \sum_{\nu = 0}^j \mathbf{E}[\bar I_i' | \bar R(j) = \nu] \bP(\bar R(j) = \nu) - \sum_{\nu = 0}^j \mathbf{E}[\bar I_{i+1} | \bar R(j) = \nu] \bP(\bar R(j) = \nu)\\
= & \sum_{\nu = 0}^j \mathbf{E}[\hat I_i' |\hat R(j) = \nu] \bP(\bar R(j) = \nu) - \sum_{\nu = 0}^j \mathbf{E}[\hat I_{i+1} |\hat  R(j) = \nu] \bP(\bar R(j) = \nu) \\ &\,\,\,\,\left(\text{\begin{tabular}{l}Given $R(j) = \nu$, reducing $t_j$ does not affect the random process after $t \geq j$,\\ due to the memoryless property.\end{tabular}}\right)\\
= & \sum_{\nu = 0}^{j-2} \mathbf{E}[\hat I_i' |\hat R(j) = \nu] \bP(\hat R(j) = \nu)+\sum_{\nu = j-1}^{j} \mathbf{E}[\hat I_i' |\hat R(j) = \nu] \bP(\bar R(j) = \nu) \\
& \,\,\,\,\,- \sum_{\nu = 0}^{j-2} \mathbf{E}[\hat I_{i+1} |\hat  R(j) = \nu] \bP(\hat R(j) = \nu)- \sum_{\nu = j-1}^{j} \mathbf{E}[\hat I_{i+1} | \hat R(j) = \nu] \bP(\bar R(j) = \nu)\\
= & \sum_{\nu = 0}^{j} \mathbf{E}[\hat I_i' |\hat R(j) = \nu] \bP(\hat R(j) = \nu)+\sum_{\nu = j-1}^{j} \mathbf{E}[\hat I_i' |\hat R(j) = \nu] (\bP(\bar R(j) = \nu)  - \bP(\hat R(j) = \nu))\\
& \,\,\,\,\,- \sum_{\nu = 0}^{j} \mathbf{E}[\hat I_{i+1} |\hat  R(j) = \nu] \bP(\hat R(j) = \nu)- \sum_{\nu = j-1}^{j} \mathbf{E}[\hat I_{i+1} | \hat R(j) = \nu] (\bP(\bar R(j) = \nu) - \bP(\hat R(j) = \nu))\\
= & \mathbf{E}[\hat I_i'] + \sum_{\nu = j-1}^{j} \mathbf{E}[\hat I_i' |\hat R(j) = \nu] (\bP(\bar R(j) = \nu) - \bP(\hat R(j) = \nu))\\
&\,\,\,\,\, - \mathbf{E}[\hat I_{i+1}] - \sum_{\nu = j-1}^{j} \mathbf{E}[\hat I_{i+1} | \hat R(j) = \nu] (\bP(\bar R(j) = \nu) - \bP(\hat R(j) = \nu))\\
= &\mathbf{E}[\hat I_i'] - \mathbf{E}[\hat I_{i+1}] + \sum_{\nu = j-1}^{j} (\mathbf{E}[\hat I_i' |\hat R(j) = \nu]- \mathbf{E}[\hat I_{i+1} | \hat R(j) = \nu]) (\bP(\bar R(j) = \nu) - \bP(\hat R(j) = \nu))\\
= &\mathbf{E}[\hat I_i] - \mathbf{E}[\hat I_{i+1}] + \sum_{\nu = j-1}^{j} (\mathbf{E}[\hat I_i' |\hat R(j) = \nu]- \mathbf{E}[\hat I_{i+1} | \hat R(j) = \nu]) (\bP(\bar R(j) = \nu) - \bP(\hat R(j) = \nu))\\
\geq & \sum_{\nu = j-1}^{j} (\mathbf{E}[\hat I_i' |\hat R(j) = \nu]- \mathbf{E}[\hat I_{i+1} | \hat R(j) = \nu]) (\bP(\bar R(j) = \nu) - \bP(\hat R(j) = \nu))\\
= & (\bP(\bar R(j) = \nu) - \bP(\hat R(j) = \nu)) \times \\
& \,\,\,\,\, ( \mathbf{E}[\hat I_i' | \hat R(j) = j] - \mathbf{E}[\hat I_{i+1}|\hat R(j) = j] - \mathbf{E}[\hat I_i'| \hat R(j) = j-1] + \mathbf{E}[\hat I_{i+1}| \hat R(j) = j-1]) . \,\,\,\text{(by (\ref{eq:nboundlmpfchangeprob}))}
\end{align*}
Now Lemma \ref{lm:nboundlm5} gives 
\[ \mathbf{E}[\hat I_i' |\hat  R(j) = j] - \mathbf{E}[\hat I_i' |\hat  R(j) = j-1] - \mathbf{E}[\hat I_{i+1} | \hat R(j) = j] + \mathbf{E}[\hat I_{i+1} |\hat  R(j) = j-1] \geq 0.\]
This proves that $\mathbf{E}[\bar I_i] \geq \mathbf{E}[\bar I_{i+1}]$. Therefore, we always have $\mathbf{E}[I_i] \geq \mathbf{E}[I_{i+1}]$ when we change $t_j$ from $j$ to $\bar t_j$ for all $j = 1,2,...,l$.
\halmos\endproof

\noindent {\bf Proof of Theorem \ref{thm:nboundlmtp}.}

\proof{Proof.}
We give a new and more detailed construction of the same set of times points as constructed in the proof of Theorem \ref{thm:nboundthm1}.

Fix any $\beta \in (0,1)$. Starting with $t_i = i$, $\forall i=1,2,...,k-1$, we run the following algorithm:
\begin{enumerate}
\item[] For $i = 0,1,...,k-2$:
\begin{enumerate}
\item If $\mathbf{E}[I_i] > 1/\beta-1$, reduce $t_{i+1}$ such that $\mathbf{E}[I_i] = 1/\beta-1$.
\item Stop if $\mathbf{E}[I_i] < 1/\beta-1$.
\end{enumerate}
\end{enumerate}

If the algorithm stops at (b) when $i=l$, we must have 
\[ \mathbf{E}[I_0] = \mathbf{E}[I_1] = \cdots = \mathbf{E}[I_{l-1}] = 1/\beta - 1\]
and, according to Lemma \ref{lm:nboundlmkey},
\begin{equation}\label{eq:nboundlmpfbbb}
1/\beta -1 = \mathbf{E}[I_{l-1}] \geq \mathbf{E}[I_l] \geq \cdots \geq \mathbf{E}[I_{k-1}].
\end{equation}

On the other hand, if the algorithm never stops at (b), we must have
\begin{equation}\label{eq:nboundlmpfaaa}
 \mathbf{E}[I_0] = \mathbf{E}[I_1] = \cdots = \mathbf{E}[I_{k-2}] = 1/\beta - 1.
 \end{equation}

If we change the value of $\beta$, the time points $t_1,t_2,...,t_{k-1}$ as the result of the algorithm must change continuously in $\beta$. This implies that $\mathbf{E}[I_{k-1}] = \int_{t_{k-1}}^k \bP(R(s) = k-1)ds$ must change continuously in $\beta$. When $\beta$ is close to $0$, we must have $\mathbf{E}[I_{k-1}] < 1/\beta -1$; when $\beta$ is close to $1$, we must have $\mathbf{E}[I_{k-1}] > 1/\beta -1$. Therefore, there must exist a $\beta$ such that, when the algorithm ends, 
\[ \mathbf{E}[I_{k-1}] = 1/\beta-1.\]
Let $\beta^*$ be such a value. 

Now the time points have met all desired conditions if (\ref{eq:nboundlmpfaaa}) holds (i.e., the algorithm never stops at (b)). If the algorithm stops at some step (b), then according to (\ref{eq:nboundlmpfbbb}),
\[ 1/\beta^* -1 = \mathbf{E}[I_{l-1}] \geq \mathbf{E}[I_l] \geq \cdots \geq \mathbf{E}[I_{k-1}] = 1/\beta^* -1\]
\[ \Longrightarrow \mathbf{E}[I_0] = \mathbf{E}[I_1] = \cdots = \mathbf{E}[I_{k-1}] = 1/\beta^* - 1,\]
which gives all desired conditions of the time points as well.
\halmos\endproof

\begin{lemma}\label{lm:nboundlm1}
\[\beta^* = \frac{1}{2} + \frac{1}{2k} \sum_{i=0}^{k-1} i\,\alpha^*_i(k).\]
\end{lemma}
\proof{Proof.}
Starting from Theorem \ref{thm:nboundthm1} we can deduce that
\[ k (1- \beta^*) = \beta^*\left[k - \sum_{i=0}^{k-1} i \bP(R(k) = i)\right] \]
\[ \Longrightarrow 2k\beta^* = k + \beta^* \sum_{i=0}^{k-1} i \bP(R(k) = i) = k + \sum_{i=0}^{k-1} i \alpha^*_i(k)\]
\[ \Longrightarrow \beta^* = \frac{1}{2} + \frac{1}{2k} \sum_{i=0}^{k-1} i \alpha^*_i(k).\]

\halmos\endproof

\noindent {\bf Proof of Lemma \ref{lm:nboundlm3}.}

\proof{Proof.}
If suffices to prove the case when $x = y+1$. We have
\begin{align*}
&\sum_{i=-l}^l P_{k-1+i-(y+1)}(\lambda) - \sum_{i=-l}^l P_{k-1+i-y}(\lambda)\\
= & P_{k-2-l-y}(\lambda) - P_{k-1+l-y}(\lambda).
\end{align*}
If $k-2-l-y < 0$, the lemma trivially holds because $P_{k-2-l-y}(\lambda) = 0$ and thus $P_{k-2-l-y}(\lambda) - P_{k-1+l-y}(\lambda) \leq 0$.

Now suppose $k - 2 - l - y \geq 0$. Then
\begin{align*}
& \frac{P_{k-2-l-y}(\lambda)}{ P_{k-1+l-y}(\lambda)}\\
= & \frac{\lambda^{ k-2-l-y}}{(k-2-l-y)!} \cdot \frac{(k-1+l-y)!}{\lambda^{k-1+l-y}}\\
= & \frac{1}{\lambda^{2l+1}} \prod_{i=-l}^l (k-1-y+i).
\end{align*}

Since $y \geq k-1-\lambda$, we must have $k-1-y\leq \lambda$ and $(k-1-y+i)(k-1-y-i)\leq \lambda^2$. This shows that $P_{k-2-l-y}(\lambda)/P_{k-1+l-y}(\lambda)\leq 1$ and thus $P_{k-2-l-y}(\lambda) - P_{k-1+l-y}(\lambda) \leq 0$.
\halmos\endproof
\begin{lemma}\label{lm:nboundlm2}
For any $l=0,1,...,k-1$, we have
\[  \sum_{i=-l}^{l} P_{k-1+i}(k) \leq \frac{1}{\beta^*}\sum_{i=0}^l \alpha_{k-1-i}^*(k).\]
\end{lemma}
\proof{Proof.}
Define
\[ R^{(i)} \equiv R(t_i) + N(k) - N(t_i), \,\,\forall i = 0,1,2,...,k.\]
Note that since the bounded process $R(t)$ is determined by $N(t)$, the random variables $R^{(i)}$'s are also determined by $N(t)$. 

Since 
\[\bP(R^{(0)} = i ) = \bP(R(t_0) + N(k) - N(t_0) = i) = \bP(N(k) = i) = P_i(k)\]
and
\[ \bP(R^{(k)} = i) = \bP(R(k) = i) = \alpha^*_i(k) / \beta^*,\]
it suffices to show that
\begin{equation}\label{eq:nboundpf1}
 \sum_{i=-l}^l \bP( R^{(j-1)} = k-1+i) \leq \sum_{i=-l}^l \bP(R^{(j)} = k-1+i)
 \end{equation}
for all $l = 0,1,...,k-1$ and $j = 1,2,...,k$.

According to the definition of the bounded process $R(t)$, if $R(t_{j-1}) + N(t_j) - N(t_{j-1}) \leq j-1$, then $R(t_j) = R(t_{j-1}) + N(t_j) - N(t_{j-1})$ and thus 
\begin{align*}
 R^{(j)} &= R(t_j) + N(k) - N(t_j) \\
 &= R(t_{j-1}) + N(t_j) - N(t_{j-1}) + N(k) - N(t_j) \\
 &= R(t_{j-1}) + N(k) - N(t_{j-1})\\
 & = R^{(j-1)}.
 \end{align*}

Therefore,
\begin{align*}
& \sum_{i=-l}^l \bP( R^{(j-1)} = k-1+i | R(t_{j-1}) + N(t_j) - N(t_{j-1}) \leq j-1) \\
= &\sum_{i=-l}^l \bP(R^{(j)} = k-1+i | R(t_{j-1}) + N(t_j) - N(t_{j-1}) \leq j-1)
\end{align*}
for all $l = 0,1,...,k-1$ and $j = 1,2,...,k$.

Now consider the case $x = R(t_{j-1}) + N(t_j) - N(t_{j-1}) > j-1$. We must have
\[ R^{(j-1)} = x + N(k) - N(t_j),\]
\[ R^{(j)} = j-1 + N(k) - N(t_j).\]
Recall that Lemma \ref{thm:nboundlmtp} gives $t_j \leq j$, so $x > j-1 \geq t_j -1 = k -1  - (k - t_j)$. We can then apply Lemma \ref{lm:nboundlm3} by further setting $y = j-1$ and $\lambda = k-t_j$ and obtain (for $x > j-1$)
\begin{align*}
& \sum_{i=-l}^l \bP( R^{(j-1)} = k-1+i | R(t_{j-1}) + N(t_j) - N(t_{j-1}) =x) \\
= &\sum_{i=-l}^l \bP(x + N(k) - N(t_j) = k-1+i | R(t_{j-1}) + N(t_j) - N(t_{j-1}) =x)\\
= &\sum_{i=-l}^l \bP(N(k) - N(t_j) = k-1+i-x | R(t_{j-1}) + N(t_j) - N(t_{j-1}) =x)\\
= & \sum_{i=-l}^l \bP(N(k-t_j) = k-1+i-x)\\
\leq &\sum_{i=-l}^l \bP(N(k-t_j) = k-1+i-(j-1))\\
= & \sum_{i=-l}^l \bP(j-1 + N(k)-N(t_j) = k-1+i | R(t_{j-1}) + N(t_j) - N(t_{j-1}) =x)\\
= & \sum_{i=-l}^l \bP( R^{(j)}= k-1+i | R(t_{j-1}) + N(t_j) - N(t_{j-1}) =x)
\end{align*}

In sum, we have shown (\ref{eq:nboundpf1}), which proves the lemma.
\halmos\endproof

\noindent {\bf Proof of Proposition \ref{prop:nboundthmy}.}
\proof{Proof.}
Combining Lemma \ref{lm:nboundlm1} and Lemma \ref{lm:nboundlm2}, we obtain
\begin{align*}
\beta^* = & \frac{1}{2} + \frac{1}{2k} \sum_{i=0}^{k-1} i\,\alpha^*_i(k)\\
= & \frac{1}{2} + \frac{1}{2k} \sum_{l=0}^{k-2}\sum_{i=0}^l \alpha^*_{k-1-i}(k)\\
\geq & \frac{1}{2} + \frac{1}{2k} \sum_{l=0}^{k-2} \beta^* \sum_{i=-l}^l P_{k-1+i}(k)\\
= & \frac{1}{2} + \frac{\beta^*}{2k} \left[ \sum_{l=0}^{k-2} \sum_{i = -l}^0 P_{k-1+i}(k) + \sum_{l=0}^{k-2} \sum_{i = 1}^l P_{k-1+i}(k) \right]\\
= & \frac{1}{2} + \frac{\beta^*}{2k} \left[\sum_{i=1}^{k-1} i P_i(k) + \sum_{i=k}^{2k-2} (2k - i-2) P_i(k) \right]\\
= & \frac{1}{2} + \frac{\beta^*}{2k} \left[\sum_{i=1}^{2k-2} i P_i(k) + \sum_{i=k}^{2k-2} (2k - 2i-2) P_i(k) \right]\\
= & \frac{1}{2} + \frac{\beta^*}{2k} \left[ \sum_{i=1}^{2k-2} i P_i(k) - 2 \sum_{i=1}^{k-1} i P_{k-1+i}(k)\right].
\end{align*}
\begin{align*}
\Longrightarrow \beta^*\geq & \frac{k}{2k - \left[ \sum_{i=1}^{2k-2} i P_i(k) - 2 \sum_{i=1}^{k-1} i P_{k-1+i}(k)\right]}\\
= &\frac{k}{k + \left[ k - \sum_{i=1}^{2k-2} i P_i(k) \right] + 2 \sum_{i=1}^{k-1} i P_{k-1+i}(k)}\\
= &\frac{k}{k + \sum_{i=2k-1}^{\infty} i P_i(k)  + 2 \sum_{i=1}^{k-1} i P_{k-1+i}(k)}\\
= &\frac{1}{1 + \frac{1}{k}\left[\sum_{i=2k-1}^{\infty} i P_i(k)  + 2 \sum_{i=1}^{k-1} i P_{k-1+i}(k)\right]}
\end{align*}
\halmos\endproof
\noindent {\bf Proof of Theorem \ref{thm:betaBound}.}
\proof{Proof.}
\begin{eqnarray*}
 \sum_{i=2k-1}^\infty i P_i(k) +2 \sum_{i=1}^{k-1} i P_{k+i-1}(k) \leq & 2 \sum_{i=1}^\infty i P_{k+i-1}(k)\\
 = & 2 \sum_{i=1}^\infty i \frac{k^{k+i-1}}{(k+i-1)!}e^{-k}\\
 = & 2 \left[ \frac{k^k e^{-k}}{(k-1)!} + \sum_{i=k}^\infty P_i(k)\right].\\
 \end{eqnarray*}
\begin{eqnarray*}
 \Longrightarrow \beta^* \geq & \frac{1}{1 + \frac{2}{k}\left[ \frac{k^k e^{-k}}{(k-1)!} + \sum_{i=k}^\infty P_i(k)\right]}\\
 = &  \frac{1}{1 + 2 \left[\frac{e^{-k} k^k}{k!}+\frac{P_{\geq k}(k)}{k} \right]}.
 \end{eqnarray*}
By Stirling's formula, 
\[ \frac{e^{-k} k^k }{k!} = \frac{1}{\sqrt{2\pi k}} + o(1/k).\]
Furthermore, since $P_{\geq k}(k) \leq 1$, we have
\[ \frac{P_{\geq k}(k)}{k} = O(1/k).\]

Thus, 
\[ \beta^* \geq \frac{1}{1 + 2 \left[\frac{e^{-k} k^k}{k!}+\frac{P_{\geq k}(k)}{k}\right]} = \frac{1}{ 1 + 2\left[ \frac{1}{\sqrt{2\pi k}} + o(1/k) + O(1/k)\right]  } = 1 - \sqrt{\frac{2}{\pi}} \frac{1}{\sqrt{k}} + O(1/k).\]

In sum, we have proved that 
\[ \frac{f_j(0,k)}{\sum_{i=1}^m x^*_{ij} r_{ij}} \geq \frac{u_0(0)}{\int_0^1 r(t)\lambda(t) dt}\geq \beta^* \geq 1 - \sqrt{\frac{2}{\pi}} \frac{1}{\sqrt{k}} + O(1/k).\]

Since $\beta^*$ increases in the capacity $k$, when $k$ is the minimum capacity of any resource, we must have
\[ \frac{\sum_{j =1}^n f_j(0,C_j)}{\sum_{j =1}^n\sum_{i=1}^m x^*_{ij} r_{ij}} \geq \min_{j \in [n]} \frac{f_j(0,C_j)}{\sum_{i=1}^m x^*_{ij} r_{ij}} \geq   \beta^* \geq 1 - \sqrt{\frac{2}{\pi}} \frac{1}{\sqrt{k}} + O(1/k).\]

This proves the competitive ratio of the Separation Algorithm.
\halmos\endproof

\noindent {\bf Proof of Theorem \ref{thm:algdomination}.}
\proof{Proof.}
Let $h(t,c)$ be the expected future reward of the Marginal Allocation Algorithm starting at time $t$ when $c = (c_1,c_2,...,c_n)$ is the vector of resources available at $t$. Let $f(t,c)=\sum_{j =1}^n f_j(t,c_j)$ be the expected future reward of the Separation Algorithm starting at time $t$ when $c$ is the vector of resources available at $t$. We will show that
\begin{equation} h(t,c) \geq f(t,c) \label{eq:twoalg}\end{equation}
for any state $(t,c)$.

Define an algorithm $\Pi^{(i)}$ as follows. For the first $i$ customers, apply the Marginal Allocation Algorithm. Afterward, for the $(i+1)$-th, $(i+2)$-th,..., customers, apply the Separation Algorithm. Let $h^{(i)}(t,c)$ be the expected future reward when algorithm $\Pi^{(i)}$ is applied starting at time $t$ with remaining inventory $c$, and \textit{assuming that no customers have arrived prior to time $t$}. We must have
\[ h^{(0)}(t,c) = f(t,c),\]
\[ \lim_{i\to \infty} h^{(i)}(t,c) = h(t,c) .\]
HJB equation for computing the expected reward of algorithm $\Pi^{(1)}$ is 
\begin{equation}\label{eq:OptimizedAlgProof2}
\frac{\partial h^{(1)}(t,c)}{\partial t} = - \sum_{i=1}^m \lambda_i(t) \left(\max\limits_{j\in [n]} \left(r_{ij} - f_j(t,c_j) + f_j(t,c_j-1)\right)^+ -  \Delta^{(1)}(t,c)\right),
\end{equation}
 where
\[ \Delta^{(1)}(t,c) \equiv h^{(1)}(t,c) - f(t,c).\]
The boundary conditions for (\ref{eq:OptimizedAlgProof2}) are $h^{(1)}(1,c) = 0$ (the length of the horizon is $1$) and $f_j(t,-1) = -\infty$, $\forall j=1,2,...,n$ and $t \in [0,1]$. Note that according to the boundary condition $f_j(t,-1) = -\infty$, the `bid price' of any resource $j$ that has no remaining inventory becomes infinity, as $f_j(t,0) - f_j(t,-1) = \infty$.

To see why (\ref{eq:OptimizedAlgProof2}) is true, consider the discrete version of (\ref{eq:OptimizedAlgProof2}). During any small period $(t, t+\delta t)$, one of the following three events must take place.
\begin{itemize}
\item No customer arrives during $(t,t+\delta t)$. Then the expected future reward $h^{(1)}(t,c)$ turns into $h^{(1)}(t+\delta t, c)$. 
\item A customer of some type $i$ arrives, but $\Pi^{(1)}$ (which applies the Marginal Allocation Algorithm to the customer) rejects the customer. We must have  
\[\max\limits_{j\in [n]} \left(r_{ij} - f_j(t + \delta t,c_j) + f_j(t + \delta t,c_j-1)\right)^+= 0,\] 
because $r_{ij}$ must be smaller than the `bid price' $f_j(t + \delta t,c_j) - f_j(t + \delta t,c_j-1)$ of all available resources $j$; we must also have that $h^{(1)}(t,c)$ turns into $f(t+\delta t,c)$, as $\Pi^{(1)}$ turns into the Separation Algorithm.
\item A customer of some type $i$ arrives and the customer is assigned to a resource $j$. In this case, the system collects reward $r_{ij}$, and $h^{(1)}(t,c)$ turns into $f(t + \delta t,c - e_j)$ as $\Pi^{(1)}$ turns into the Separation Algorithm, where $e_j$ is the unit vector with the $j$-th position being $1$. Note that 
\[f(t + \delta t, c - e_j) = f(t + \delta t, c) + f_j(t + \delta t, c_j-1) - f_j(t + \delta t, c_j). \]
\\
Then mathematically, we can combine the second and the third bullet points, and say that when a customer of type $i$ arrives, the expected future reward $h^{(1)}(t,c)$ turns into total \emph{current and future} reward
\[ f(t + \delta t, c) + \max\limits_{j\in [n]} \left(r_{ij} - f_j(t + \delta t,c_j) + f_j(t + \delta t,c_j-1)\right)^+.\]
\end{itemize}
In sum, the recursive equation for $h^{(1)}(t,c)$ is
\begin{align*}
 h^{(1)}(t, c) = & (1 - \sum_{i=1}^m \lambda_i(t)\delta t) h^{(1)}(t + \delta t, c) \\
 &+ \sum_{i=1}^m \lambda_i(t) \delta t\left(  f(t + \delta t, c) + \max\limits_{j\in [n]} \left(r_{ij} - f_j(t + \delta t,c_j) + f_j(t + \delta t,c_j-1)\right)^+ \right).
 \end{align*}
 Letting $\delta t \to 0$ we obtain (\ref{eq:OptimizedAlgProof2}).

Therefore,
\begin{equation}\frac{\partial h^{(1)}(t,c)}{\partial t} \leq  \frac{\partial f(t,c)}{\partial t} + \sum_{i=1}^m \lambda_i(t)  \Delta^{(1)}(t,c).  \end{equation}
This equation implies that, if at some time $t_0$ we have $\Delta^{(1)}(t_0,c)< 0$ or equivalently
\begin{equation} \label{eq:OptimizedAlgProofb}
h^{(1)}(t_0,c) - f(t_0,c) < 0,
 \end{equation}
then we must have
\begin{equation}  \frac{\partial h^{(1)}(t,c)}{\partial t} <  \frac{\partial f(t,c)}{\partial t}, \,\,\, \forall t \in (t_0,1] \end{equation}
and
\begin{equation} \label{eq:OptimizedAlgProofa}
 h^{(1)}(t,c) < f(t,c), \,\,\, \forall t \in (t_0,1].
 \end{equation}
However, since we know that $h^{(1)}(1,c) = f(1,c) = 0$, (\ref{eq:OptimizedAlgProofa}) cannot be true, and thus (\ref{eq:OptimizedAlgProofb}) cannot be true. Therefore, we have proved 

\begin{equation}\label{eq:OptimizedAlgProof4}
h^{(1)}(t,c) \geq f(t,c), \,\,\, \forall t \in [0,1].
 \end{equation}
 
Next, we show that 
\begin{equation}\label{eq:OptimizedAlgProof5}
h^{(i)}(t,c) \geq h^{(i-1)}(t,c), \,\,\, \forall t\in [0,1]
\end{equation}
 by induction on $i$. 
 
Equation (\ref{eq:OptimizedAlgProof4}) already proves the base case $i=1$. Suppose for some $\bar i > 1$, (\ref{eq:OptimizedAlgProof5}) holds for all $i < \bar i$. Now we show that it also holds for $i = \bar i$. By definition, for any $\bar i > 1$, algorithms $\Pi^{(\bar i)}$ and $\Pi^{(\bar i-1)}$ must allocate the first customer in the same way. Thus, $\Pi^{(\bar i)}$ and $\Pi^{(\bar i-1)}$ earn the same reward from the first customer, and then transit into the same state. After that first customer,  $\Pi^{(\bar i)}$ continues to apply $\Pi^{(\bar i-1)}$ pretending that no customer has ever arrived, while $\Pi^{(\bar i-1)}$ continues to apply $\Pi^{(\bar i-2)}$. By induction, the expected future reward of  $\Pi^{(\bar i-1)}$ is at least that of $\Pi^{(\bar i-2)}$. Therefore, the expected future reward of $\Pi^{(\bar i)}$ is at least that of $\Pi^{(\bar i-1)}$.

Thus, we have proved (\ref{eq:OptimizedAlgProof5}). It immediately follows that 
\[ h^{(\infty)}(t,c) \geq h^{(0)}(t,c).\]
\halmos
\endproof

\bibliographystyle{ormsv080}
\bibliography{myrefs}{}

\end{document}